\documentclass[dvipsnames]{amsart}

\usepackage{amsmath,amssymb,amsthm}

\usepackage{hyperref}
\usepackage{units}
\usepackage{pgfplots}
\usepackage{fullpage}
\usepackage{todonotes}
\usepackage{relsize}
\usepackage{rotating}
\usepackage{lscape}
\usepackage{mathtools}
\usepackage{lipsum}
\usepackage{graphicx}
\usepackage{afterpage}
\newtheorem{theorem}{Theorem}[section]
\newtheorem{lemma}[theorem]{Lemma}	
\newtheorem{proposition}[theorem]{Proposition}

\theoremstyle{definition}
\newtheorem{definition}[theorem]{Definition} 
\newtheorem{remark}[theorem]{Remark}	
\newtheorem{example}[theorem]{Example}

\usepackage{mathptmx}       
\usepackage{helvet}         
\usepackage{courier}        
\usepackage{type1cm}        
\usepackage{graphicx}        
\usepackage{multicol}        
\usepackage[bottom]{footmisc}

\newcommand{\Rn}{\mathbb R^n}

\newcommand{\rn}[1]{{\mathbb R}^{#1}}
\newcommand{\R}{\mathbb R}

\newcommand{\G}{\mathbb G}

\newcommand{\he}[1]{{\mathbb H}^{#1}}

\newcommand{\cov}[1]{{\bigwedge\nolimits^{#1}{\mfrak g}}}
\newcommand{\vet}[1]{{\bigwedge\nolimits_{#1}{\mfrak g}}}

\newcommand{\covw}[2]{{\bigwedge\nolimits^{#1,#2}{\mfrak g}}}

\newcommand{\vetguno}[1]{{\bigwedge\nolimits_{#1}{\mfrak g_1}}}

\newcommand{\scal}[2]{\langle {#1} , {#2}\rangle}

\newcommand{\Scal}[2]{\langle {#1} \vert {#2}\rangle}
\newcommand{\scalp}[3]{\langle {#1} , {#2}\rangle_{#3}}

\newcommand{\ccheck}{{\vphantom i}^{\mathrm v}\!\,}

\newcommand{\mc}{\mathcal }

\newcommand{\mfrak}{\mathfrak}

\newcommand{\gf}{\varphi}

\newcommand{\gt}{\theta}

\newcommand{\N}{\mathbb N}

\title[Comparison of three possible Laplacians]{Comparing three possible hypoelliptic Laplacians on the 5-dimensional Cartan group via div-Curl type estimates} 
 
\author[A. Baldi and F. Tripaldi]{Annalisa Baldi and Francesca Tripaldi}
\address[A. Baldi]%
{Universit\`a di Bologna, Dipartimento
di Matematica, Piazza di
Porta S.~Donato 5, 40126 Bologna, Italy. } 
\email{annalisa.baldi2@unibo.it}

\address[F. Tripaldi]%
{Centro De Giorgi,
SNS,
Piazza dei Cavalieri 3,  56 126 Pisa, Italy.} 
\email{francesca.tripaldi@sns.it}

\date{}

\begin{document}

\maketitle

\begin{abstract}
On general Carnot groups, the definition of a possible hypoelliptic Hodge-Laplacian on forms using the Rumin complex has been considered in \cite{rumin_gafa,rumin_cras}, where the author introduced a 0-order pseudodifferential operator on forms.  However, for questions regarding regularity for example, where one needs sharp estimates, this 0-order operator is not suitable. Up to now, there have only been very few attempts to define hypoelliptic Hodge-Laplacians on forms that would allow for such sharp estimates. Indeed, this question is rather difficult to address in full generality, the main issue being that the Rumin exterior differential $d_c$ is not homogeneous on arbitrary Carnot groups. 

In this note, we consider the specific example of the free Carnot group of step 3 with 2 generators, and we introduce three possible definitions of hypoelliptic Hodge-Laplacians. We compare how these three possible Laplacians can be used to obtain sharp div-curl type inequalities akin to those considered by Bourgain \& Brezis and Lanzani \& Stein for the de Rham complex, or their subelliptic counterparts obtained by Baldi, Franchi \& Pansu for the Rumin complex on Heisenberg groups.
\end{abstract}

\section{Introduction}




Div-curl inequalities in $\rn n$ read as follows 
\begin{eqnarray}
\label{0a} &\|u\|_{L^{n/(n-1)}(\rn n)}\le C\, \big( \|du \|_{L^1(\rn n)} + \|\delta u \|_{L^1(\rn n)}\big) 
\end{eqnarray}
when $u$ is a  compactly supported
smooth differential $h$-form with $h\neq 1,n-1$. Here, $d$ is the exterior differential, and $\delta$ the exterior
co-differential. 
In addition, denoting by $\mc H^1(\R^n)$ the Hardy space, one also has
\begin{align*}
    &\|u\|_{L^{n/(n-1)}(\rn n)}\le C\, \big( \|du \|_{L^1(\rn n)} + \|\delta u \|_{\mc H^1(\R^ n)}\big) \ \text{ when }h=1\ \text{ and }\\
    &\|u\|_{L^{n/(n-1)}(\rn n)}\le C\, \big( \|du \|_{\mc H^1(\rn n)} + \|\delta u \|_{L^1(\rn n)}\big)\ \text{ when }h=n-1\,.
\end{align*}
These estimates have been proved by 
 Lanzani \& Stein \cite{LS} and by Burgain \& Brezis
\cite{BB2007}, and are applied to the study
of div-curl systems and of more general Hodge systems in $\Rn$.
 They contain
in particular the famous Burgain-Brezis inequality \cite{BB2003}, \cite{BB2004} for divergence-free vector fields in $\rn n$  (see also 
\cite{vS2004}).

We stress that the previous inequality for $h=1$ (and analogously for $h=n-1$) fails if we replace the Hardy norm with
the $L^1$-norm, i.e. an estimate like \eqref{0a} is {\it false} for $1$-forms (and $(n-1)$-forms), with a counterexample given by E.M. Stein in \cite{stein_diff}, p. 191.

Obtaining analogous estimates in the setting of Carnot groups has proven much more difficult. In \cite{CvS2009},  Chanillo \& Van Schaftingen extended these Burgain-Brezis inequalities to a class of vector fields
in Carnot groups and, 
outside the Euclidean setting, estimates akin to \eqref{0a} have been proved in \cite{BF7} and \cite{BFP1} in the setting of the Heisenberg groups for differential forms belonging to the so-called Rumin complex.

The Rumin complex  $(E_0^*,d_c)$ is a subcomplex of de Rham complex $(\Omega^*,d)$ that M. Rumin constructed on Carnot groups \cite{rumin_grenoble,rumin_cras,rumin_jdg}. It is a complex, homotopic to the de Rham complex, that better fits the geometry of the group. Recently, Rumin's theory has been fruitfully used to address several questions in differential
geometry, as well as in pde's theory in Carnot groups (and more generally in sub-Riemannian structures). In the particular case of Heisenberg groups, with homogeneous dimension $Q=2n+2$, one can obtain the following sharp estimates if $u$ is a form in $E_0^h$,  $1<h<2n$ and $h\neq n,n+1$, so that
\begin{align}\label{0b} 
&\|u\|_{L^{Q/(Q-1)}(\he{n}, E_0^h)}\le   C\big( \|d_cu\|_{L^1(\he{n}, E_0^{h+1})} + \|\delta_cu \|_{L^1(\he{n}, E_0^{h-1})}\big)\,.
\end{align}
 
On the other hand, if $h=n,n+1$, then the differential $d_c$, and the co-differential $\delta_c$ respectively, have order 2. The estimates \eqref{0b} then take the form 
\begin{align*} 
&\|u\|_{L^{Q/(Q-2)}(\he{n}, E_0^n)}\le   C\big( \|d_cu \|_{L^1(\he{n}, E_0^{n+1})} + \|d_c\delta_cu \|_{L^1(\he{n}, E_0^{n})}\big) ;
\\
&\|u\|_{L^{Q/(Q-2)}(\he{n}, E_0^{n+1})}\le   C\big(\|\delta_c d_cu \|_{L^1(\he{n}, E_0^{n+1})} + \|\delta_cu \|_{L^1(\he{n}, E_0^{n})}  \big) \,.
\end{align*}
For a more precise statement, including the cases when $h=1,2n$, see Theorem 1.3 of \cite{BFP1}. The proof of this result relies on the precise estimates for the fundamental solution of the hypoelliptic Laplace-type operators defined for forms of any degree for Heisenberg groups. These Laplacians have been introduced by Rumin in \cite{rumin_jdg}.

 In this work, we focus on extending this type of limiting Sobolev inequalities in the spirit of \eqref{0a} to the Cartan group. This is achieved by defining three different classes of hypoelliptic Hodge-Laplacian operators on the Rumin complex, all of which turn out to be suitable for obtaining such estimates. More precisely, we aim to compare the possible div-curl type estimates that we obtain for each class of Laplacians.
 
 The Cartan group is the 5-dimensional free nilpotent Carnot group whose Lie algebra is spanned by $(X_1,X_2,X_3=[X_1,X_2],X_4=[X_1,X_3],X_5=[X_2,X_3])$. 
This group, and more generally 5-dimensional nilmanifolds with generic\footnote{We recall that a 2-dimensional distribution is said to be generic, or $(2,3,5)$ on a nilmanifold $M$, if it admits the Cartan group as a local frame.} rank two distribution, have been extensively studied as they arise in several contexts. For example, the mechanical system of a surface rolling without slipping and twisting on another surface gives rise to a 5-dimensional configuration space equipped with a rank two distribution encoding the no slipping and twisting condition. This provides a basic example of $(2,3,5)$ nilmanifolds (under suitable conditions see e.g. Section 4.4 in \cite{BryantHsu}).  The Cartan group also provides an example of conformal geometry which is naturally associated with ordinary differential equations \cite{CapNeusser,Nurowski}.

So far, there has been extensive work dealing with hypoelliptic Hodge-Laplace operators of order 12 on the Cartan group  \cite{Haller2022,haller2023regularized,haller2023analytic}. In these quoted works, these operators (referred to as Rumin-Seshadri operators) are constructed using the Rumin differential\footnote{More precisely, their construction of the subcomplex relies on the BGG machinery, which has been recently shown to coincide with the Rumin complex on homogeneous groups, and hence on Carnot groups \cite{F+T1}}, with application to heat kernel expansion and analytic torsion, as a generalisation of the work done by Rumin on contact manifolds.

In our paper, these homogeneous order 12 Laplacians on forms will be denoted by $\Delta_\G$, since they are essentially (up to a possible rescaling of the metric) the hypoelliptic Hodge-Laplace operators in  \cite{Haller2022,haller2023regularized,haller2023analytic}. In addition, we also define two other different families of Hodge-Laplace operators. They are again constructed using Rumin differentials. A first family will be denoted by $\Delta_{A}$ which turn out of order 6, except for degrees 0 and 5, where they have order 2 (see Definition \ref{rumin laplacian}). The second family denoted by $\Delta_{R}$ differs from $\Delta_A$ in degrees $2$ and $3$, where $\Delta_R$ has differential order 12 (see Definition \ref{RS Laplacian}). Moreover, $\Delta_R$ also differs from the Laplacian $\Delta_\G$ which has homogeneous order 12 in each degree. 

All the operators $\Delta_\G$, ${\Delta}_{R}$ and $\Delta_{A}$ are hypoelliptic, as shown in Proposition \ref{hypo 2}. Moreover, the operators $\Delta_{A}$ have differential orders 6 and 2 strictly less than the homogeneous dimension of the Cartan group (which is 10). Hence, $\Delta_A$ has a fundamental solution to which it is possible to apply the theory of singular integrals of Folland and Stein \cite{folland_stein}, as shown in Theorem 3.1 of \cite{BFT3}. For the other two Laplacians, the differential order can be 12, higher than the homogenous dimension $Q=10$. It is however still possible to obtain estimates of the derivative of the fundamental solutions of both $\Delta_\G$ and $\Delta_{R,h}$ ($h=2,3$) thanks to a recent result obtained by Van Schaftingen and Yung for vector-valued operators (see Theorem 3.3 in \cite{VanSchaftingenYung}).

Rumin's theory needs a quite technical introduction that is sketched
in Section \ref{rumin forms} to make the paper self-constained. 
The main properties of $(E_0^*,d_c)$ 
can be summarised in the following
points:
\begin{itemize}
\item Rumin $1$-forms are horizontal $1$--forms, i.e. forms that are dual
 of horizontal vector fields.
 
   \item The  ``intrinsic'' exterior differential $d_c$ on a smooth function is
   its horizontal differential (that is dual operator of the gradient along a basis
   of the horizontal bundle). Instead, when acting on forms of degree $\ge 1$, the differential $d_c$ can be identified with a matrix-valued operator where each entry is a \textit{homogeneous left-invariant differential operator of order} $\ge 1$. We stress that, in general, different entries have different differential orders, and so $d_c$, viewed as a matrix-valued operator,  \textit{is not homogeneous}.
      \item The complex $(E_0^*,d_c)$ is exact and self-dual under Hodge $\ast$-duality.
   \end{itemize}

Among all Carnot groups, Heisenberg groups $\he n$  provide the simplest examples of noncommutative
Carnot groups (of step 2). Recently, Rumin's Laplacians in Heisenberg groups, or in sub-Riemannian contact manifolds with bounded geometry have been used to get Poincar\'e inequalities and sharp estimates on $\ell^{q,p}$-cohomology. Indeed, thanks to its scale
invariance, Rumin's Laplacian allows to apply the theory of singular integral operators, 
 to get the sharp exponent in Heisenberg groups, also in degree $h = n, n+1$ where the operator $d_c$ is of order two (see \cite{BFP3,PT,BFP2,BFP5,BFP6}. Moreover, Sobolev-Gaffney inequalities associated to the Rumin complex in this same setting have also relied on such hypoelliptic Laplacians (see \cite{BF7,BFP1,BTTr}).

Unfortunately, the geometry of arbitrary Carnot groups is radically different to that of  Heisenberg groups. 
One of the difficulties that one encounters is related to the fact that
 the structure of intrinsic differential forms is much more complex and, in general, the exterior differential is not homogeneous as a matrix-valued operator. As a consequence, it is not obvious how to construct a hypoelliptic Laplacian on forms (as instead is the case in the Heisenberg groups). This
lack of homogeneity is connected to the notion of {\it weight} of a differential form (see Definition \ref{peso peso}). As far as we know, except for the 0-order Hodge-Laplacians introduced by Rumin e.g. in \cite{rumin_cras,rumin_grenoble}, the only attempt to construct hypoelliptic Hodge-Laplacians in the case where the Rumin differentials are \textit{not homogeneous} can be found in \cite{BFTr}, where we discuss the case of the Engel group (which is the model of a class of so-called filiform groups) and the $7$-dimensional quaternionic Heisenberg group. 

However, in the case of the Cartan group, the Rumin differentials $d_c$ are homogeneous on each degree. This property has been used in \cite{DaveHallerBGG} to construct the Rumin-Seshadri operators $\Delta_\G$, and it's also the reason behind the homogeneity of the Hodge-Laplacians $\Delta_R$. In our paper, we again use the homogeneity of $d_c$ to construct a brand new Laplace-type hypoelliptic operator on forms $\Delta_A$. All three operators $\Delta_\G$, $\Delta_R$, and $\Delta_A$ can be used to obtain sharp estimates of div-curl type on the Cartan group. These results are contained in Theorems \ref{H2}, \ref{C2}, and \ref{H2cor}.

The fundamental solution of the hypoelliptic $\Delta_A$ plays a crucial role in obtaining the estimates of Theorem \ref{H2}. On the other hand, the different div-curl type estimates of Theorem \ref{C2} follow from the definition of Laplacians $\Delta_\G$ and $\Delta_R$ given in Definitions \ref{DH Laplacian} and \ref{RS Laplacian}, respectively. Even though for arbitrary forms the estimates of Theorem \ref{C2} are not equivalent to those of Theorem \ref{H2} (obtained using $\Delta_A$), it is important to stress that if one instead considers a differential form $u\in E_0^h$ which is either closed or co-closed, the relative estimates from all three Laplacians coincide. For example, in the case of a co-closed form, the statement of Theorem \ref{H2cor} can be roughly summarised as follows.

\medskip

\noindent\textit{Denote by $(E_0^*,d_c)$ the complex of Rumin forms in the 5-dimensional Cartan group $\G$.
Then there exists $C>0$ such that, for any compactly supported smooth Rumin $h$-form
$u$ so that $\delta_cu=0$,
we have}
 \begin{align*}
&\|u\|_{L^{Q/(Q-1)}(\G)\hphantom{, E_0^5}}\le   C\|d_cu \|_{L^1(\G, E_0^{1})} \quad & \mbox{if $h=0$;}
\\
&\|u\|_{L^{Q/(Q-3)}(\G, E_0^1)}\le  C \|d_cu \|_{L^1(\G,E_0^2)}  \quad & \mbox{if $h=1$;}
\\
&\|u\|_{L^{Q/(Q-2)}(\G, E_0^2)}\le   C \| d_cu \|_{L^1(\G, E_0^3)}  \quad& \mbox{if $h=2$.}
\end{align*}
\textit{where $Q=10$ is the homogeneous dimension of $\G$. The estimates for the other degrees can be readily obtained by Hodge duality.
}

Moreover, interestingly, when considering interpolation-type $L^{p}$-spaces, the relative div-curl estimates read as follows (see Section \ref{ultima} for the precise statement, here we write only the case $h=1$ and $h=2$):
\medskip

\noindent{\it  Denote by $(E_0^*,d_c)$ the complex of intrinsic forms in $\G$.
Then there exists $C>0$ such that for any $h$-form
$u\in \mc D(\G, E_0^h)$, $h=1,2$, such that
$$
\left\{\begin{aligned}
 d_c u = f \\ \delta_c u = g&
\end{aligned}
\right.
$$
we have
 \begin{align*}
&\|u\|_{L^{Q/(Q-3)}+L^{Q/(Q-1)}(\G, E_0^1)}\le  C\big( \|f \|_{L^1(\G,E_0^2)} + \|g \|_{\mc H^1(\G)}\big) \quad & \mbox{if $h=1$;}
\\
&\|u\|_{L^{Q/(Q-3)}+L^{Q/(Q-2)}(\G, E_0^2)}\le   C\big( \| f \|_{L^1(\G, E_0^3)} + \|g \|_{L^1(\G, E_0^1)}\big) \quad& \mbox{if $h=2$.}
\end{align*}
}
\medskip

The paper is organised as follows. In Section \ref{notation section}, we give some basic definitions related to Carnot groups and multilinear algebra and present the main analytic results that we will use in later sections. Section \ref{rumin forms} contains a quick review of the Rumin complex. Moreover, in Section \ref{Rumin conti espliciti}, we address the particular example of the 5-dimensional Cartan group $\G$ and present some explicit computations of the Rumin complex in this specific setting. In addition, in the same section, we write the definitions for the three families of Hodge-Laplace operators $\Delta_\G$, $\Delta_{R}$, and $\Delta_A$. We provide an explicit proof of the hypoellipticity of $\Delta_A$, and we recall some properties of its fundamental solution. Finally, in Section \ref{section 5}, we present the statements of our limiting Sobolev-type estimates (see Theorems \ref{H2}, \ref{C2} and \ref{H2cor}). These results follow from the properties of the fundamental solutions of the three associated Laplace-type operators, together with a result of Chanillo and Van Schaftingen, stated in Theorem \ref{chanillo_van5}, that can be applied after an appropriate rephrasing of the closeness of forms in terms of symmetric tensors (see Subsection \ref{rephrasing}).
Finally, in Subsection \ref{ultima}, we provide an alternative family of $L^1$-estimates using interpolation-type $L^p$-spaces.

\section{Notations and preliminary results on Carnot group}\label{notation section}
A Lie group $G$ is a smooth manifold endowed with the smooth mappings
\begin{align*}
    G\times G\ni (x,y)\mapsto xy\in G\ \text{ and }\ G\ni x\mapsto x^{-1}\in G
\end{align*}
satisfying, for all $x,y,z\in G$, the properties
\begin{itemize}
    \item [1.] $x(yz)=(xy)z$;
    \item [2.] $ex=xe=x$;
    \item [3.] $xx^{-1}=x^{-1}x=e$,
\end{itemize}
where $e\in G$ is the unit element of the group $G$.

A {\it Carnot group $\G$ of
step $\kappa$}  is a connected, simply connected
Lie group whose Lie algebra
${\mathfrak{g}}$ admits a {\it step $\kappa$ stratification}, i.e.
there exist linear subspaces $V_1,...,V_\kappa$ such that
\begin{equation}\label{stratificazione}
{\mathfrak{g}}=V_1\oplus...\oplus V_\kappa,\quad [V_1,V_i]=V_{i+1},\quad
V_\kappa\neq\{0\},\quad V_i=\{0\}{\textrm{ if }} i>\kappa,
\end{equation}
where $[V_1,V_i]$ is the subspace of ${\mathfrak{g}}$ generated by
the commutators $[X,Y]$ with $X\in V_1$ and $Y\in V_i$. Let
$m_i=\dim(V_i)$, for $i=1,\dots,\kappa$ and $h_i=m_1+\dots +m_i$ with
$h_0=0$ and, clearly, $h_\kappa=n$. Choose a basis $e_1,\dots,e_n$ of
$\mathfrak{g}$ adapted to the stratification, that is such that
$$e_{h_{j-1}+1},\dots,e_{h_j}\;\text {is a basis of}\; V_j\;\text{
for each}\; j=1,\dots, k.$$
Let $X_1,\dots,X_{n}$ be the family
of left-invariant vector fields such that
$X_i(0)=e_i$. Given (\ref{stratificazione}), the subset
$X_1,\dots,X_{m_1}$
generates by commutations all the other vector fields; we will
refer to $X_1,\dots,X_{m_1}$ as {\it generating  vector fields} of
the group. The exponential map is a one-to-one map from $\mathfrak
g$ onto $\G$, i.e. any $p\in\G$ can be written in a unique way as
$p=\exp(p_1X_1+\dots+p_nX_n)$. Using these {\it exponential
coordinates}, we identify $p$ with the n-tuple $(p_1,\dots,p_n)\in
\R^n$ and we identify $\G$ with $(\Rn,\cdot)$ where the explicit
expression of the group operation $\cdot$ is determined by the
Campbell-Hausdorff formula (see \cite{folland,BLU}).

The subbundle of the tangent bundle $T\G$ that is spanned by the
vector fields $X_1,\dots,X_{m_1}$ plays a particularly important
role in the theory, it is called the {\it horizontal bundle}
$H\G$. The sections of $H\G$ are called {\it horizontal
sections}, and a vector of $H_x\G$ is an {\it horizontal vector}. 

A subriemannian
structure is defined on $\G$, endowing each fiber of $H\G$ with a
scalar product $\scal{\cdot}{\cdot}_{x}$ and with a norm
$\vert\cdot\vert_x$ that make the basis $X_1(x),\dots,X_{m_1}(x)$
an orthonormal basis. 

Adopting the standard notation, given $f\in\mathcal{D}(\G)$, we denote by $\nabla_\G f$ the usual horizontal gradient of the function $f$, i.e.
\begin{align*}
    \nabla_\G f=(X_1f,X_2f,\ldots,X_{m_1}f)\,.
\end{align*}


Two important families of automorphisms of $\G$ are the 
group translations and the group dilations of $\G$. For
any $x\in\G$, the {\it (left) translation} $\tau_x:\G\to\G$ is defined
as $$ z\mapsto\tau_x z:=x\cdot z. $$ For any $\lambda >0$, the
{\it dilation} $\delta_\lambda:\G\to\G$, is defined as
\begin{equation}\label{dilatazioni}
\delta_\lambda(x_1,...,x_n)=
(\lambda^{d_1}x_1,...,\lambda^{d_n}x_n),
\end{equation} where $d_i\in\N$ is called {\it homogeneity of
the variable} $x_i$ in
$\G$ (see \cite{folland_stein} Chapter 1) and is defined as
\begin{equation}\label{omogeneita2}
d_j=i \quad\text {whenever}\; h_{i-1}+1\leq j\leq h_{i},
\end{equation}
hence $1=d_1=...=d_{m_1}<
d_{{m_1}+1}=2\leq...\leq d_n=\kappa.$
In particular, the generating vector fields are homogeneous of
degree $1$ with respect to group dilations.

Several distances equivalent to $d_c$ have been used in the literature. 
If $p=(\tilde p_1,\dots,
\tilde p_\kappa)\in \R^{m_1}\times\cdots\times\R^{m_\kappa}=\Rn$, then we set
\begin{equation}\label{distanzainfinito}
\rho(p):=\max\{\varepsilon_j||\tilde{p}_j||_{\R^{m_j}}^{1/j}\, ,
j=1,\dots,\kappa\}. \end{equation} Here $\varepsilon_1=1$, and
$\varepsilon_{2},\dots\varepsilon_{\kappa} \in (0,1)$ are suitable
positive constants depending on the group structure (see Theorem
5.1 in \cite{FSSC_JGA}). We can consider the following 
gauge distance  $$d(x,y):=\rho(y^{-1}\cdot x)\,.$$

\noindent The metric $d$
 well behaves with respect to left
translations and dilations, that is
\begin{equation*}
d(z\cdot x,z\cdot y)=d(x,y)\quad ,\quad
d(\delta_\lambda(x),\delta_\lambda(y))=\lambda d(x,y)
\end{equation*}
for
$x,y,z\in\G$ and $\lambda>0$.


\bigskip

The integer
\begin{equation}\label{dimensione}
Q=\sum_{i=1}^\kappa i\,\text{dim}\,V_i
\end{equation}
is the
{\it homogeneous dimension} of $\G$. It is also the
Hausdorff dimension of $\Rn$ with respect to the  $d$.

The $n$-dimensional Lebesgue measure $\mathcal L^n$, is the Haar
measure of the group $\G$. Hence if $E\subset\Rn$ is measurable,
then $ \mathcal L^n(x\cdot E)= \mathcal L^n(E)$ for all $x\in\G$.
Moreover, if $\lambda>0$ then $\mathcal
 L^n(\delta_\lambda(E))=\lambda^Q \mathcal L^n(E)$. 
All the spaces $L^p(\G)$ that appear throughout this paper are
defined with respect to the $\mathcal L^n$ Lebesgue measure.

Following e.g. \cite{folland_stein}, we can define a group
convolution in $\G$: if, for instance, $f\in\mc D(\G)$ and
$g\in L^1_{\mathrm{loc}}(\G)$, we set
\begin{equation}\label{group convolution}
f\ast g(p):=\int f(q)g(q^{-1}\cdot p)\,dq\quad\mbox{for $q\in \G$}.
\end{equation}
We remind that, if (say) $g$ is a smooth function and $L$
is a left-invariant differential operator, then
$
L(f\ast g)= f\ast Lg.
$
We remind also that the convolution is again well defined
when $f,g\in\mc D'(\G)$, provided at least one of them
has compact support (as customary, we denote by
$\mc E'(\G)$ the class of compactly supported distributions
in $\G$ identified with $\rn {n}$, the dual of $\mc E(\G)$,
the space of smooth functions). 

We now recall the notion of {\it kernel of order $\mu$}.
Following \cite{folland}, a kernel of order $\mu$ is a 
homogeneous distribution of degree $\mu-Q$
(with respect to group dilations $\delta_r$ as in (\ref{dilatazioni}), see
\cite{folland_stein}),
that is smooth outside of the origin. Alternatively, we shall also write \textit{kernel of type $\mu$} to mean the same thing.

\begin{proposition}\label{kernel}
Let $K\in\mc D'(\Omega)$ be a kernel of order $\mu$. If $\mu>0$, then $K\in L^1_{\mathrm{loc}}(\G)$.
Moreover, $X_\ell K$ is a kernel of order $\mu-1$ for
any horizontal derivative $X_\ell K$.
\end{proposition}

In Proposition 1.11 in \cite{folland}, the following result is proved.

\begin{theorem}\label{Folland 75 thm}
    Suppose $0<\mu <Q$, $1<p<Q/\mu$ and $\frac{1}{q}=\frac{1}{p}-\frac{\mu}{Q}$. Let $K$ be a kernel of type $\mu$. If $u\in L^p(\G)$ the convolutions $u\ast K$ and $K\ast u$ exist a.e. and are in $L^q(\G)$ and there is a constant $C_p>0$ such that 
	$$
	\|u\ast K\|_q\le C_p\|u\|_p\quad \mathrm{and}\quad \| K\ast u\|_q\le C_p\|u\|_p\,
	$$
	
\end{theorem}

Following \cite{folland_stein}, we also adopt the following
multi-index notation for higher-order derivatives. If $I =
(i_1,\dots,i_{n})$ is a multi--index, we set 
\begin{equation}\label{centrata}
X^I=X_1^{i_1}\cdots
X_{n}^{i_{n}}.
\end{equation}
 By the Poincar\'e--Birkhoff--Witt
theorem (see, e.g. \cite{bourbaki}, I.2.7), the differential
operators $X^I$ form a basis for the algebra of left invariant
differential operators in $\G$. Furthermore, we set
$|I|:=i_1+\cdots +i_{n}$ the order of the differential
operator $X^I$, and   $d(I):=d_1i_1+\cdots +d_n i_{n}$ its
degree of homogeneity with respect to group dilations.
 From the Poincar\'e--Birkhoff--Witt theorem, it follows, in particular, that any
homogeneous linear differential operator in the horizontal
derivatives can be expressed as a linear combination of the
operators $X^I$ of the special form above. Thus, we can
restrict ourselves to consider operators of the special form
$X^I$.

 We finish this section recalling the definition of free Carnot group.

\begin{definition}\label{free} Let $m_1 \ge 2$ and $ \kappa \ge 1$  be fixed integers.
We say that $\mathfrak{ f}_{m_1,\kappa} $ is the free Lie algebra with $m_1$ generators $X_1,\dots, X_{m_1}$ and
nilpotent of step $\kappa$ if:
\begin{itemize}
\item[i)] $\mathfrak f_{m_1,\kappa}$ is a Lie algebra generated by its elements $X_1,\dots, X_{m_1}$;
\item[ii)]  $\mathfrak f_{m_1,\kappa}$ is nilpotent of step $\kappa$;
\item[iii)]  for every Lie algebra $\mathfrak n$ nilpotent of step $\kappa$ and for every map $\phi$ from the set
$\{X_1,\dots, X_{m_1}\}$ to $\mathfrak n$, there exists a (unique) homomorphism of Lie algebras
$\Phi$
from
$\mathfrak f_{m_1,\kappa}$ to $\mathfrak n$ which extends $\phi$.
\end{itemize}

A Carnot group  $\G$ is said free if its Lie algebra $\mathfrak g$ is isomorphic to 
a free Lie algebra.
\end{definition}
The Cartan group is the free Carnot group with $m_1=2$ and $\kappa=3$.
\begin{example}[The Cartan group]  \label{Cartan group def}
    Let us consider now the free group of step 3 with 2
 generators, also known as the Cartan group, i.e. 
the Carnot group whose Lie algebra
is
$\mathfrak g=V_1\oplus V_2\oplus V_3$ with
$V_1=\mathrm{span}\;\{X_1,X_2\}$, $V_2=\mathrm{span}\;\{X_3\}$, and
$V_3=\mathrm{span}\;\{X_4,X_5\}$, the only non zero commutation relations
being
$$
[X_1,X_2]=X_3\quad , \quad [X_1,X_3]=X_4\quad , \quad [X_2,X_3]=X_5.
$$
In exponential coordinates, the group can be identified with
$\rn{5}$, and an explicit representation of the
vector fields is
\begin{eqnarray*}
X_1=\partial_1 &,&
X_2=\partial_2+x_1\partial_3+\frac{x_1^2}{2}
\partial_4+x_1x_2\partial_5\\
X_3=\partial_3+x_1\partial_4+x_2\partial_5 &,&
X_4=\partial_4\quad , \quad X_5=\partial_5.
\end{eqnarray*}
For an alternative expression of $X_i$, we refer to group $N_{5,2,3}$ in \cite{cornucopia}.

\end{example}

\subsection{Multilinear algebra}
The dual space of $\mfrak g$ is denoted by $\cov 1$ and indicate by $\scalp{\cdot}{\cdot}{}$ also the inner product in $\cov 1$. Given $\lbrace X_1,\ldots,X_n\rbrace$ an orthonormal basis of $\mathfrak g$, its dual basis is denoted by 
$\{\theta_1,\cdots, \theta_{n}\}$. 

Following Federer (see \cite {federer} 1.3), the exterior algebras of 
$\mfrak g$ and of $\cov 1$ are the graded algebras indicated as
$\displaystyle\vet \ast =\bigoplus_{h=0}^{n} \vet h
$
 and 
 $\displaystyle  \cov \ast 
=\bigoplus_{h=0}^{n} \cov h
$
  where
$
 \vet 0 = \cov 0 =\R
$
and, for $1\leq h \leq n$,
\begin{equation*}
\begin{split}
         \vet h& :=\mathrm {span}\{ X_{i_1}\wedge\dots \wedge X_{i_h}: 1\leq
i_1< \dots< i_h\leq n\},   \\
         \cov h& :=\mathrm {span}\{ \gt_{i_1}\wedge\dots \wedge \gt_{i_h}:
1\leq i_1< \dots< i_h\leq n\}.
\end{split}
\end{equation*}
The elements of $\vet h$ and $\cov h$ are called  \emph{$h$-vectors} and \emph{
$h$-covectors}. 

\medskip
We denote by $\Theta^h$ the basis $\{ \gt_{i_1}\wedge\dots \wedge \gt_{i_h}
:
1\leq i_1< \dots< i_h\leq n\}$
of $  \cov h$.
We remind that
$
\dim \cov h = \dim \vet h = \binom{n}{h}.
$

The dual space $\bigwedge^1(\vet h)$ of 
$\vet h$ can be naturally identified with $\cov h$.
The action of a $h$-covector $\gf$ on a $h$-vector $v$ is denoted as $\Scal
\gf v$.

The inner product $\scalp{\cdot}{\cdot}{} $ extends canonically to 
$\vet h $ and to $\cov h$ making the bases
$X_{i_1}\wedge\dots \wedge X_{i_h}$ and $\gt_{i_1}\wedge\dots \wedge
\gt_{i_h}$ orthonormal.

Starting from $\vet *$ and
$\cov *$, by left translation, we can define now two families of vector bundles over $\G$ (still
denoted by $\vet *$ and
$\cov *$) (see \cite{BFTT} for details). Smooth sections of these
vector bundles are said respectively vector fields and differential forms.

\begin{definition} If $0\le h\le n$ and $1\le p\le \infty$,
we denote by $L^{p}(\G,\cov{h})$
the space of all sections of $\cov{h}$ such that their
components with respect to the basis $\Theta^h$ belong to
$L^{p}(\G)$. Clearly, this definition
is independent of  the choice of the basis itself.
\end{definition}

We conclude this section by recalling the notion of
Hodge duality, see  \cite{federer} 1.7.8.

For $0\le h\le n$, we denote by
$\star : \cov h \longleftrightarrow \cov{n-h}$,
the $\star$-Hodge isomorphism associated to the scalar product $\scalp{\cdot}{\cdot}{}$ and the volume form $\theta_1\wedge\cdots\wedge\theta_n$.

We recall that, $\star\star \varphi=(-1)^{(n-h)h}\varphi$ for any $\varphi\in \cov h$.
\subsection{Weight of forms}\label{subsection weights of forms}
The notion of weight of a differential form plays a key role when considering the de Rham complex on a Carnot group.

\begin{definition}\label{peso peso} If $\alpha\in \cov 1$, $\alpha\neq 0$, 
 we say that $\alpha$ has \emph{pure weight $p$}, and we write
$w(\alpha)=p$, if $\alpha^\natural\in V_p$. 
More generally, if
$\alpha\in \cov h$, we say that $\alpha$ has pure weight $p$ if $\alpha$ is
a linear combination of covectors $\theta_{i_1}\wedge\cdots\wedge\theta_{i_h}$
with $w(\theta_{i_1})+\cdots + w(\theta_{ i_h})=p$.

In particular, the volume form has weight $Q$ (the homogeneous
dimension of the group).
\end{definition}
\begin{remark}[Remark 2.4 in \cite{BFTT}]\label{orthogonality}
If $\alpha,\beta \in \cov h$ and $w(\alpha)\neq w(\beta)$, then
$\scal{\alpha}{\beta}=0$. 

\end{remark}
From this Remark, we readily have the orthogonal decomposition $
\cov h = \bigoplus_{p=M_h^{\mathrm{min}}}^{M_h^{\mathrm{max}}} \covw {h}{p}$,
where $\covw {h}{p}$ is the linear span of the $h$--covectors of weight $p$
and $M_h^{\mathrm{min}}$, $M_h^{\mathrm{max}}$ are respectively the smallest
and the largest weight of left-invariant $h$-covectors.

Since the elements of the basis $\Theta^h$ have pure weights, a basis of
$ \covw {h}{p}$ is given by $$\Theta^{h,p}:=\Theta^h\cap \covw {h}{p}\,.$$
In other words, the basis $\Theta^h =\cup_p \Theta^{h,p}$ is a basis adapted to
the filtration of $ \cov h$.

We denote by  $\Omega^{h,p} $ the vector space of all
smooth $h$--forms in $\G$ of pure weight $p$, i.e. the space of all
smooth sections of $\covw {h}{p}$. We have
\begin{equation}\label{deco forms}
\Omega^h = \bigoplus_{p=M_h^{\mathrm{min}}}^{M_h^{\mathrm{max}}} \Omega^{h,p}.
\end{equation}

\begin{lemma}[Section 2.1 in \cite{rumin_grenoble}]\label{d0 left} We have
$d(\covw{h}{p})\subset \covw{h+1}{p}$, i.e., if $\alpha\in \covw{h}{p}$ is a left invariant 
$h$-form of weight p with $d\alpha \neq 0$,
then $w(d\alpha)=w(\alpha)$.
\end{lemma}
\begin{definition}\label{di}
Let now $\alpha = \sum_{\theta^h_i\in\Theta^{h,p}}\alpha_{i}\, \theta_i^h
\in \Omega^{h,p}$ be a (say) smooth form
of pure weight $p$. Then
we can write
$$
d\alpha=d_0\alpha+d_1\alpha+\dots+d_{\kappa}\alpha,
$$
where $$d_0\alpha = \sum_{\theta^h_i\in\Theta^{h,p}}\alpha_id\theta_i^h$$ does not increase the weight,
$$d_1\alpha= \sum_{\theta^h_i\in\Theta^{h,p}}\sum_{j=1}^{m_1}(X_j\alpha_i)\theta_j\wedge\theta_i^h$$
increases the weight of $1$, and, more generally,
$$d_p\alpha= \sum_{\theta^h_i\in\Theta^{h,p}}\sum_{X_j\in V_p}(X_j\alpha_i)\theta_j\wedge\theta_i^h,$$
when $p=0,1,\dots,\kappa$.
In particular, $d_0$ is an algebraic operator.

\end{definition}

\section{The Rumin complex on Carnot groups}\label{rumin forms}
The notion of intrinsic form in Carnot groups is due to M. Rumin
(\cite{rumin_grenoble}, \cite{rumin_cras}). A more extended presentation of the results of this section
can be found in \cite{BFTT}, \cite{FT4}.

\begin{definition}[M. Rumin]\label{E0}
If $0\le h\le n$ we set
$$
E_0^h:=\ker d_0\cap\ker \delta_0 = \ker d_0\cap (\mathcal{R} (d_0))^{\perp}
\subset \Omega^h
$$
Here $\delta_0:\cov{h+1}\to \cov{h}$ is the adjoint of $d_0$ with
respect to our fixed scalar product.

In the sequel, we refer to the elements of $E_0^h$ as to {\it Rumin or intrinsic
$h$-forms on $\G$}.
\end{definition}

Since the construction of $E_0^h$ is left invariant, this space of forms can be seen as the space of sections of a fiber subbundle of $\bigwedge^h\mathfrak g$, generated by left translations and still denoted by $E_0^h$. In particular, $E_0^h$ inherits from $\bigwedge^h\mathfrak g$ the scalar product on the fibers.

As a consequence, we can obtain a left invariant orthonormal basis
$\Xi_0^h=\{\xi_j\}$ of $E_0^h$ such that
\begin{equation*}
\Xi^h_0 = \bigcup_{p=M_h^{\mathrm{min}}}^{M_h^{\mathrm{max}}} \Xi^{h,p}_0,
\end{equation*}
where $ \Xi^{h,p}_0:= \Xi^h\cap\covw{h}{p}$ is a left-invariant orthonormal basis of
$E_0^{h,p}$. All the elements of $\Xi^{h,p}_0$ have pure
weight $p$.
Without loss of generality, the indices $j$ of $\Xi_0^h=\{\xi_j^h\}$
are ordered once and for all in an increasing way according
to the weight of the respective element of the basis.
 
Correspondingly, the set of indices $\{1,2,\dots,\dim E^h_0\}$ can be written as
the union of finite (possibly empty) sets of indices
$$
\{1,2,\dots,\dim E^h_0\} = \bigcup_{p= M_h^{\mathrm{min}}}^{M_h^{\mathrm{max}}}
I^h_{p},
$$
where
$$
j\in I^h_{0,p}\quad\mbox{if and only if}\quad
\xi^h_j\in \Xi^{h,p}_0.
$$
 
Without loss of generality, if $m:=\dim V_1$, we can take
$$
\Xi_0^1 = \Xi_0^{1,1}  = \{\theta_1,\ldots,\theta_m\}.
$$
 
Once the basis $\Xi_0^h$ is chosen, the space $\mathcal{D}(\G,E_0^h)$,  can be identified with
$\mc D(\G)^{\dim E_0^h}$.

\medskip
In the last few years, there have been several
 instances that show the
``naturalness'' of using the Rumin complex instead of de Rham complex for Carnot groups. For example, we can quote \cite{FT4}, Theorem 3.16, where the authors study the naturalness of $d_c$ in terms of homogeneous homomorphisms of the group $\G$, or \cite{BF4} where the authors show that $d_c$ appears, in the spirit of the Riemannian approximation, as limit. 

\begin{definition} If $0\le h\le n$ and $1\le p\le \infty$,
we denote by $L^p(\G,E_0^h)$
the space of all sections of $E_0^h$ such that their
components with respect to the basis $\Xi^h_0$ belong to
$L^p(\G)$, endowed with its natural norm. Clearly, this definition
is independent of  the choice of the basis itself.
\end{definition}

\begin{lemma}[\cite{BFTT}, Lemma 2.11]\label{d_0}
If $\beta\in\cov{h+1}$, then there exists a unique $\alpha\in
\cov h\cap (\ker d_0)^\perp$ such that
$$
\delta_0 d_0\alpha = \delta_0\beta.\quad\mbox{We set}\quad
\alpha :=d_0^{-1}\beta.
$$
 In particular
$$
\alpha =d_0^{-1}\beta\quad\mbox{if and only if}\quad
d_0\alpha-\beta\in \mc R(d_0)^\perp.
$$
Moreover
\begin{itemize}
\item[i)] $(\ker d_0)^\perp= \mc R(d_0^{-1})$;
\item[ii)] $d_0^{-1}d_0 = Id$ on $(\ker d_0)^\perp$;
\item[iii)] $d_0d_0^{-1} -Id: \cov {h+1}\to \mc R(d_0)^\perp$.
\end{itemize}
\end{lemma}

The following  theorem summarizes the construction of 
the intrinsic differential $d_c$ (for details, see \cite{rumin_grenoble}
and \cite{BFTT}, Section 2) .

\begin{theorem}\label{main rumin new}
The de Rham complex $(\Omega^\ast,d)$ splits in the direct sum of two sub-complexes, denoted by Rumin as $(E^\ast,d)$ and $(F^\ast,d)$. Moreover, we have:
\begin{itemize}
\item[i)] Let $\Pi_E$ be the projection on $E$ along $F$, then
for any $\alpha\in E_0^{h,p}$, if we denote by $(\Pi_E\alpha)_j$
the component of $\Pi_E\alpha$ of weight $j$, then
\begin{equation}\label{rumin 6}\begin{split}
(\Pi_E\alpha)_p &=\alpha\\
(\Pi_E\alpha)_{p+k+1} &= -d_0^{-1}\big(
\sum_{1\le\ell\le k+1}d_\ell (\Pi_E\alpha)_{p+k+1-\ell}\big).
\end{split}\end{equation}
Notice that $\alpha\to (\Pi_E\alpha)_{p+k+1}$ is an homogeneous
differential operator of order $k+1$ in the horizontal derivatives.
\item[ii)] $\Pi_{E}$ is a chain map, i.e.
$$
d\Pi_{E} = \Pi_{E}d.
$$
\end{itemize}

Set now
 $$d_c=\Pi_{E_0}\, d\,\Pi_{E}: E_0^h\to E_0^{h+1}, \quad h=0,\dots,n-1\,,$$
 where $\Pi_{E_0}$ is the orthogonal projection of $\Omega^\ast$ onto $E_0^\ast$.

 Then $(E_0^\ast,d_c)$ is a complex computing the de Rham cohomology of the underlying Carnot group such that
\begin{itemize}
\item[iii)] the differential $d_c$
		acting on $h$-forms can be identified, with respect to the
		bases $\Xi_0^h$ and $\Xi_0^{h+1}$, with a matrix-valued
		differential operator $L^h:= \big(L^h_{i,j}\big)$. 
		If $j\in I^h_{p}$ and $i\in I^{h+1}_{q}$, then the
		$L^h_{i,j}$'s are homogeneous left invariant differential
		operator of order $q-p\ge 1$ in the horizontal derivatives,
		and 
		$L^h_{i,j}=0$  if $j\in I^h_{0,p}$ and $i\in I^{h+1}_{0,q}$, with
		$q-p<1$
\item[iv)] the space of Rumin forms $E_0^\ast$ is closed under Hodge-star duality and the $L^2$-formal adjoint $\delta_c$ of $d_c$ on $E_0^h$ satisfies 
	\begin{equation}\label{jan 9 eq:2}
		\delta_c = (-1)^{n(h+1)+1} \star d_c\star. 
  \end{equation}
\end{itemize}
\end{theorem}

\section{The Rumin complex for the 5-dimensional Cartan group and Hypoelliptic Laplace operators on the group}\label{Rumin conti espliciti}
From now on, we denote by $\G$ the free group of step 3 with 2
 generators, also known as the Cartan group, introduced in Example \ref{Cartan group def}.
 
Following Example B.7 of \cite{BFTT}, denote by $\theta_1,\dots,\theta_5$ the dual left invariant forms. 
An orthonormal basis
of $E_0^1$ is given by 
\begin{align}\label{base 1 forme di Rumin}
    \Xi_0^1=\Xi_0^{1,1}=\lbrace\theta_1,\theta_2\rbrace\,.
\end{align}

We have $d\theta_1=d\theta_2=0$ and
\begin{align}\label{d0 cartan}
    d\theta_3=d_0\theta_3=-\theta_1\wedge\theta_2,\quad
d\theta_4=d_0\theta_4=-\theta_1\wedge\theta_3,\quad
d\theta_5=d_0\theta_5=-\theta_2\wedge\theta_3.
\end{align}
We explicitly note that, according to Definition \ref{peso peso}, we have the following weights:
\begin{align*}
    w(\theta_1)=w(\theta_2)=1\,,\,w(\theta_3)=2\,,\,w(\theta_4)=w(\theta_5)=3\,.
\end{align*}
The volume form $\theta_1\wedge\theta_2\wedge\theta_3\wedge\theta_4\wedge\theta_5$ is also denoted by $dV$.

Using the notations introduced in Subsection \ref{subsection weights of forms},
 we can
chose an orthonormal basis of $\cov h$, $h=1,\ldots,5$ as follows:

\begin{itemize}
\item[$\mathbf{h=1}$:] $\Theta^{1,1} = \{\theta_1,\theta_2\}
$, $\Theta^{1,2} = \{\theta_3\}$, and $\Theta^{1,3} = \{\theta_4,\theta_5\}$.
\item[$\mathbf{h=2}$:] $\Theta^{2,2} =
\{\theta_1\wedge \theta_2\}$,
$\Theta^{2,3}= \{\theta_1\wedge \theta_3, \theta_2\wedge \theta_3\}$, 
$\Theta^{2,4}= \{\theta_1\wedge \theta_4,  \theta_1\wedge \theta_5,\theta_2\wedge \theta_4,\theta_2\wedge \theta_5\}$,
$\Theta^{2,5}=\{\theta_3\wedge \theta_4,  \theta_3\wedge \theta_5\}$,
$\Theta^{2,6}= \{ \theta_4\wedge \theta_5\}$.
 \item[$\mathbf{h=3}$:] 
 $\Theta^{3,4} =
\{\theta_1\wedge \theta_2\wedge\theta_3\}$,
 $\Theta^{3,5} = 
\{\theta_1\wedge \theta_2\wedge\theta_4, 
\theta_1\wedge \theta_2\wedge\theta_5\}$,
$\Theta^{3,6} =
\{\theta_1\wedge \theta_3\wedge\theta_4, 
\theta_1\wedge \theta_3\wedge\theta_5,
\theta_2\wedge \theta_3\wedge\theta_4,
\theta_2\wedge \theta_3\wedge\theta_5\}$,
$\Theta^{3,7} = 
\{\theta_1\wedge \theta_4\wedge\theta_5,
\theta_2\wedge \theta_4\wedge\theta_5\}$,
$\Theta^{3,8} = 
\{\theta_3\wedge \theta_4\wedge\theta_5\}$.
\item[$\mathbf{h=4}:$] $\Theta^{4,7}=\lbrace \theta_1\wedge\theta_2\wedge\theta_3\wedge\theta_4,\theta_1\wedge\theta_2\wedge\theta_3\wedge\theta_5\rbrace$, $\Theta^{4,8}=\lbrace\theta_1\wedge\theta_2\wedge\theta_4\wedge\theta_5\rbrace$, $\Theta^{4,9}=\lbrace\theta_1\wedge\theta_3\wedge\theta_4\wedge\theta_5,\theta_2\wedge\theta_3\wedge\theta_4\wedge\theta_5\rbrace$.
\item[$\mathbf{h=5}:$] $\Theta^{5,10}=\lbrace\theta_1\wedge\theta_2\wedge\theta_3\wedge\theta_4\wedge\theta_5\rbrace$. 
\end{itemize}
One can readily obtain that
\begin{align*}
    \Xi_0^1=\Xi_0^{1,1}=\lbrace\theta_1,\theta_2\rbrace
\end{align*}
and by Hodge duality
\begin{align*}
    \Xi_0^4=\Xi_0^{4,9}=\lbrace\theta_1\wedge\theta_3\wedge\theta_4\wedge\theta_5,\theta_2\wedge\theta_3\wedge\theta_4\wedge\theta_5\rbrace\,.
\end{align*}

Moreover, as computed in Example B.7 of \cite{BFTT}, orthonormal basis for the space of intrinsic forms of $E_0^2$ is given by
\begin{align}\label{base 2 forme di Rumin}
    \Xi_0^2=\Xi_0^{2,4}=\{\theta_1\wedge\theta_4,\frac{1}{\sqrt{2}}(\theta_2\wedge\theta_4+\theta_1\wedge\theta_5),\theta_2\wedge\theta_5\}.
\end{align}
By Hodge duality, an orthonormal basis of $E_0^3$ is given by
\begin{align}\label{base 3 forme di Rumin}
    \Xi_0^3=\Xi_0^{3,6}=\lbrace\theta_1\wedge\theta_3\wedge\theta_4,\frac{1}{\sqrt{2}}(\theta_1\wedge\theta_3\wedge\theta_5+\theta_2\wedge\theta_3\wedge\theta_4),\theta_2\wedge\theta_3\wedge\theta_5\rbrace\,.
\end{align}
\begin{remark}\label{Im d0 sulle  forme} Let us now present the action of the map $d_0\colon\Omega^h\to\Omega^{h+1}$ according to the weights of the covectors.
\begin{itemize}
    \item [i)]$h=1$: the map $d_0$ acting on 1-forms $d_0^{(1)}\colon \Omega^1\to\Omega^2$ acts in such a way that $\mathcal{R}(d_0^{(1)})=\Omega^{2,2}\oplus\Omega^{2,3}$;
    \item [ii)] $h=2$: the map $d_0^{(2)}\colon \Omega^2\to\Omega^3$ is such that $\mathcal{R}(d_0^{(2)})=span_{C^\infty(\G)}\lbrace d_0^{(2)}(\theta_1\wedge\theta_5-\theta_2\wedge\theta_4)\rbrace\oplus \Omega^{3,5}\oplus\Omega^{3,6}=span_{C^\infty(\G)}\lbrace\theta_1\wedge\theta_2\wedge\theta_3\rbrace\oplus\Omega^{3,5}\oplus\Omega^{3,6}$
    \item [iii)] $h=3$: the map $d_0^{(3)}\colon \Omega^3\to\Omega^4$ is such that $\mathcal{R}(d_0^{(3)})=\Omega^{4,7}\oplus\Omega^{4,8}$.
\end{itemize}

\end{remark}

For general $h=0,\ldots, 4$, we should denote by $d_c^{(h)} $ the differential operators $d_c^{(h)}:E_0^{h}\to E_0^{h+1}$. To avoid cumbersome notation, when clear from the context, we shall remove the superscript $(h)$ indicating the degree of forms it acts on. 

We recall that the differentials $d_c$ are left invariant and 
homogeneous with respect to the group dilations. The differential $d_c^{(h)}$ is a first order homogeneous operator
in the horizontal derivatives in degree $h=0$ and $h=4$, it is a third
	order homogeneous operator in degree $h=1$ and $h=3$, and it is a second order homogeneous operator in degree $h=2$. 
As remarked in Theorem \ref{main rumin new}-iii), once we fix the bases for the spaces of forms $E_0^h$, it is possible to express each operator $d_c\colon E_0^h\to E_0^{h+1}$ in matrix form $(L_{i,j})$. We will now consider the ordered bases $\Xi_0^h$ listed above.
We have:

\begin{itemize}
	\item $d_{c}: {E}_0^0\longrightarrow{E}_0^1$ can be seen in matrix form as
	\begin{align*}
		d_{c}=\begin{pmatrix}
		X_1\\X_2
	\end{pmatrix}
	\end{align*}	
	\item
	$d_{c}:{E}_0^1\longrightarrow{E}_0^2$ can be expressed as
	\begin{align*}
		d_{c}=\begin{pmatrix}
		&-X_1^2X_2-X_1X_3-X_4\  & X_1^3  \\ &-\sqrt{2}\big(X_1X_2^2+X_5)\ & \sqrt{2}(X_2X_1^2-X_4) \\ &-X_2^3\ & X_2^2X_1-X_2X_3-X_5
	\end{pmatrix}
	\end{align*}
	
	\item
	$d_{c}:{E}_0^2\longrightarrow{E}_0^3$ is given by
	\begin{align*}
		d_{c}=\begin{pmatrix}
		-X_1X_2-X_3 & \frac{1}{\sqrt{2}}X_1^2 & 0\\ -\frac{1}{\sqrt{2}}X_2^2 & -\frac{3}{2}X_3 & \frac{1}{\sqrt{2}}X_1^2\\
		0 & -\frac{1}{\sqrt{2}}X_2^2 & X_2X_1-X_3
	\end{pmatrix}
 	\end{align*}
	
	\item
	$d_c:{E}_0^3\longrightarrow{E}_0^4$ can be expressed as
	\begin{align*}
		d_c=\begin{pmatrix}
		    {\left(X_{1} X_{2} + X_{3}\right)} X_{2} - X_{5}  & \sqrt{2}(-X_1^2X_2+X_4) & X_{1}^{3}\\
      X_{2}^{3}  & -\sqrt{2} (X_2^2X_1+X_5) & X_2X_{1}^2 -X_3X_1+ X_{4}
		\end{pmatrix}
	\end{align*}
 \item $d_c\colon E_0^4\longrightarrow E_0^5$ can be expressed as
 \begin{align*}
     d_c=\begin{pmatrix}
         -X_{2} & X_{1}
     \end{pmatrix}
 \end{align*}
\end{itemize}

An idea of how to get the explicit expression for the operators $d_c\colon E_0^h\to E_0^{h+1}$ can be found in the proof of Theorem \ref{pierre}.

By Hodge duality, one can find the explicit expressions for the codifferentials $\delta_c$. Once expressed in terms of the ordered bases $\Xi_0^h$, the matrix form of $\delta_c$ can be expressed as the transpose of the matrix $(L_{i,j})$ of differentials $d_c$ (up to a sign). We recall the action of the Hodge-$\star$ operator is expressed in matrix forms as
\begin{itemize}
    \item $\star_1\colon E_0^1\to E_0^4$ and $\star_4\colon E_0^4\to E_0^1$ have the form
    \begin{align*}
        \star_1=\begin{pmatrix}
            0 & -1\\ 1 & 0
        \end{pmatrix}\ \text{ and }\ \star_4=\begin{pmatrix}
            0 & 1\\ -1 & 0
        \end{pmatrix}\,,
    \end{align*}
    \item $\star_2\colon E_0^2\to E_0^3$ and $\star_3\colon E_0^3\to E_0^2$ have the form
    \begin{align*}
        \star_2=\begin{pmatrix}   
 0 & 0 &  1  \\
 0 & -1 &  0  \\
 1 & 0 &  0     \end{pmatrix}\ \text{ and }\ \star_3=\begin{pmatrix}   
 0 & 0 &  1  \\
 0 & -1 &  0  \\
 1 & 0 &  0     \end{pmatrix}
    \end{align*}
\end{itemize}
Therefore, using the formula \eqref{jan 9 eq:2}, the codifferential is expressed as $$\delta_c=(-1)^{5h}\star_{5-h+1} d_c^{(5-h)}\star_h\,\colon E_0^h\to E_0^{h-1}$$ so that
\begin{itemize}
	\item $\delta_{c}:{E}_0^1\longrightarrow{E}_0^0$ has the form
	\begin{align*}
		\delta_{c}=\begin{pmatrix}
		    -X_1 & -X_2
		\end{pmatrix}
	\end{align*}
	
	\item
	$\delta_{c}:{E}_0^2\longrightarrow{E}_0^1$ has the form:
	\begin{align*}
		\delta_{c}=\left(\begin{array}{ccc}
 {\left(X_{2} X_{1} - X_{3}\right)} X_{1} + X_{4} &  \sqrt{2}(X_2^2X_1+X_5) & X_{2}^{3}  \\
 -X_{1}^{3} &  \sqrt{2}(-X_1^2X_2+ \, X_{4}) & -{\left(X_{1} X_{2} + X_{3}\right)} X_{2} + X_{5}  
\end{array}\right)
	\end{align*}
	
	\item
	$\delta_{c}:{E}_0^3\longrightarrow{E}_0^2$ has the form:
	\begin{align*}
		\delta_{c}=\left(\begin{array}{ccc}
 -X_{2} X_{1} + X_{3} &  -\frac{1}{\sqrt{2}} \, X_{2}^{2} & 0  \\
 \frac{1}{\sqrt{2}} \, X_{1}^{2}  & \frac{3}{2} \, X_{3} & -\frac{1}{\sqrt{2}} \, X_{2}^{2}  \\
 0 & \frac{1}{\sqrt{2}} \, X_{1}^{2}  & X_{1} X_{2} + X_{3} 
\end{array}\right)
	\end{align*}
 \item $\delta_c\colon E_0^4\to E_0^3$ has the form:
 \begin{align*}
     \delta_c=\begin{pmatrix}
         -X_2^2X_1+X_2X_3+X_5 & -X_2^3\\
         \sqrt{2}(X_2X_1^2-X_4) & \sqrt{2}(X_1X_2^2+X_5)\\
         -X_1^3 & -X_1^2X_2-X_1X_3-X_4
     \end{pmatrix}
 \end{align*}
\item $\delta_c\colon E_0^5\to E_0^4$ has the form:
\begin{align*}
    \delta_c=\begin{pmatrix}
        X_2\\-X_1
    \end{pmatrix}\,.
\end{align*}
	
\end{itemize}
\subsection{Definition of three possible Laplacians following the Rumin-Seshadri approach} In this subsection, we propose three possible constructions of hypoelliptic Laplace-type operators on forms on the Cartan group. The three definitions are inspired by the paper \cite{Rumin-Seshadri}, where Rumin and Seshadri define a Laplacian of fourth order for differential forms of any degree on Heisenberberg groups in order to define a subRiemannian version of the classical analytic torsion. An analogous construction to the one in \cite{Rumin-Seshadri} has also been proposed by Dave and Haller in \cite{DaveHallerBGG} for the particular case of Carnot groups where the $d_c^{(h)}$ are homogeneous for all $h$. In both of the papers already quoted, the Hodge-Laplacians are homogeneous left-invariant differential operators on forms of the same order across all degrees. In particular, the ``Rumin-Seshadri'' operator considered in \cite{DaveHallerBGG} has always order 12 on the Cartan group.
Here we recall the definition of these 12-order Laplacians applied to our context:
\begin{definition}\label{DH Laplacian}
\begin{align*}
    \Delta_{\G,h}=
		\left\{
		\begin{array}{lcl} 
(\delta_c d_c)^6\quad &\mbox{if } & h=0;\\
(d_c\delta_c)^6+(\delta_cd_c)^2\quad&\mbox{if }&h=1;\\ (d_c\delta_c )^2+(\delta_cd_c)^3\quad& \mbox{if } & h=2;\\
			(d_c\delta_c)^3+(\delta_cd_c )^2\quad&\mbox{if }& h=3;
			\\(\delta_cd_c)^6+(d_c\delta_c)^2 \quad &\mbox{if }  & h=4;\\ (d_c\delta_c)^6
    \quad &\mbox{if }  & h=5\,.
		\end{array}
		\right.
\end{align*}
Notice that $-\Delta_{\G,0}$ has order $12$ and is NOT the usual sub-Laplacian on $\G$, i.e. $\Delta_0:=X_1^2+X_2^2$. Also for forms of higher degree, the Laplacians always have differential order 12.
\end{definition}
As proved in \cite{DaveHallerBGG}, for any $h=0,\cdots, 5$ the operators  $\Delta_{\G,h}$ are hypoelliptic self-adjoint.

\medskip

In this paper we propose, in addition, two other possible classes of Laplacians. First, the operators $\Delta_R$ below have different differential order depending on the degree of the forms they act on, and they are defined exclusively via the Rumin differentials $d_c$ and co-differentials $\delta_c$, in the spirit of the initial works of Rumin to define the Hodge-Laplacian on contact manifolds \cite{rumin_jdg}. In the second case, we define Laplace-type operators $\Delta_A$ of lower differential order, and the order stays the same across all degrees (except on functions and top-degree forms, see below for details). 

\begin{definition}\label{RS Laplacian}
Here, we consider the following class of Laplacians for the Cartan group:
\begin{align*}
    \Delta_{R,h}=
		\left\{
		\begin{array}{lcl} 
(d_c\delta_c)^3+\delta_c d_c\quad &\mbox{if } & h=0,1;			\\ (d_c\delta_c )^2+(\delta_cd_c)^3\quad& \mbox{if } & h=2;\\
			(d_c\delta_c)^3+(\delta_cd_c )^2\quad&\mbox{if }& h=3;
			\\d_c\delta_c+(\delta_c d_c)^3 \quad &\mbox{if }  & h=4,5\,.
		\end{array}
		\right.
\end{align*}
Notice that now $-\Delta_{R,0}=\Delta_0=X_1^2+X_2^2$ is the usual sub-Laplacian on $\G$. Moreover, $\Delta_{R,1}$ and $\Delta_{R,4}$ have order 6 as differential operators, while on 2 and 3-forms the Laplacians have differential order 12.
\end{definition}

We now introduce the last class of Laplace operators. 
In particular, we stress the new part of the definition below deals with the operators acting on forms of degree 2 and 3. We will denote such operators by $\Delta_{A,h}$ and we stress again that they have the same differential order across all degrees $h$ (except for $h=0,5$). To achieve this, we introduce new operators acting on 2 and 3-forms which, after composing with the $d_c$ and $\delta_c$, produce Laplacians of homogeneous differential order 6. Apart from degrees $h=2,3$, $\Delta_R$ and $\Delta_A$ coincide. 

\begin{definition}\label{rumin laplacian} 
	In $\G$, it is possible to define the following homogeneous Laplace-type operators $\Delta_{A,h}$ acting on $E_0^h$,  by setting
	\begin{equation*}
		\Delta_{A,h}=
		\left\{
		\begin{array}{lcl} 
			(d_c\delta_c)^3+\delta_c d_c\quad &\mbox{if } & h=0,1;
			\\ d_c\delta_c +\delta_cA_{\Delta}d_c\quad& \mbox{if } & h=2;\\
			d_cA_\Delta\delta_c+\delta_cd_c \quad&\mbox{if }& h=3;
			\\d_c\delta_c+(\delta_c d_c)^3 \quad &\mbox{if }  & h=4,5\,.
		\end{array}
		\right.
	\end{equation*}
	where $$A_\Delta:=-\Delta_{0}\,I_3$$ with $I_3$ is the $3\times 3$ identity matrix and $\Delta_{0} = X_1^2+X_2^2$ is again the usual subLaplacian of
$\G$. 
	
\end{definition}

Once a basis  of $E_0^h$
is fixed, the operator $\Delta_{\G,h}$ can be identified with a matrix-valued map, still denoted
by $\Delta_{\G,h}$
\begin{equation}\label{matrix form}
	\Delta_{\G{},h} = (\Delta_{\G{},h}^{ij})_{i,j=1,\dots,N_h}: \mc D'(\G, \rn{N_h})\to \mc D'(\G, \rn{N_h}),
\end{equation}
where $\mc D'(\G, \rn{N_h})$ is the space of vector-valued distributions on $\G$, and $N_h$ is
the dimension of $E_0^h$ (see \cite{BFT3}). The same is true for $\Delta_{A,h}$ and $\Delta_{R,h}$.
\begin{remark}
    By construction, it is straightforward to see that also the operators $\Delta_{R,h}$ and $\Delta_{A,h}$ are self-adjoint.

The same is true for $\Delta_A$, since 
 $$\Delta_{A,3}=\star_2 \Delta_{A,2}\star_3\,.$$

Indeed, from the definition of $\Delta_{A,2}$, we get that
 \begin{align*}
     \star \Delta_{A,2}\star=\star d_c \delta_c\star+\star\delta_cA_\Delta d_c\star=\star d_c\star\star \delta_c\star+\star\delta_c\star\star A_\Delta\star\star d_c\star=\delta_c d_c+d_c\star A_\Delta\star\delta_c
 \end{align*}
 and using the matrix form of $\star_2\colon E_0^2\to E_0^3$ and $\star_3\colon E_0^3\to E_0^2$, we get
 \begin{align*}
    \star_3 A_\Delta\star_2=&\begin{pmatrix}   
 0 & 0 &  1  \\
 0 & -1 &  0  \\
 1 & 0 &  0     \end{pmatrix} \begin{pmatrix}
     -\Delta_{0}& 0 & 0\\
     0 &-\Delta_{0}&  0\\
     0 & 0 &-\Delta_{0}
 \end{pmatrix}\begin{pmatrix}   
 0 & 0 &  1  \\
 0 & -1 &  0  \\
 1 & 0 &  0     \end{pmatrix}=A_\Delta
 \end{align*}
\end{remark}

\begin{proposition}\label{hypo 2}The operators $\Delta_{\G,h}$, $\Delta_{R,h}$ and $\Delta_{A,h}$ are hypoelliptic for each $h=0,\ldots,5$.
\end{proposition}

The proof of the hypoellipticity of the operators $\Delta_{R,h}$ and $\Delta_{\G,h}$ is essentially covered by Lemma 2.14 of \cite{DaveHallerBGG}. The hypoellipticity of the new operators $\Delta_{A,h}$ when $h=2,3$ is shown in the next section (see Proposition \ref{hypo}).

\medskip

Since for any $h=0,1,\cdots,5$ , the operators $\Delta_{A, h}$ are homogeneous hypoelliptic operators of order $6<Q=10$, we can apply the following result.

{ \begin{theorem}[see \cite{BFT3}, Theorem 3.1]\label{global solution}
		There exist
		\begin{equation}\label{numero}
			K_{ij}\in\mc D'(\G\cap \mc
			C^\infty(\G \setminus\{0\})
			\qquad\mbox{for $i,j =1,\dots,N_h$,}
		\end{equation} where, again, $N_h := \mathrm{dim}\, E_0^h$,
		with the following properties:
		\begin{enumerate}
			\item[i)]  the $K_{ij}$'s are
			kernels of type $6$,
			for
			$i,j
			=1,\dots, N_h$.
			
			\item[ii)] when $\alpha = (\alpha_1,\dots, \alpha_{N_h}) \in
			\mc D(\G,\rn {N_h})$ ,
			if we set
			\begin{equation}\label{numero2}
				\Delta_{A,h}^{-1} \alpha:=
				\big(    
				\sum_{j}\alpha_j\ast  K_{1j},\dots,
				\sum_{j}\alpha_j\ast  K_{N_hj}\big),
			\end{equation}
			then 
			$$ 
			\Delta_{A{},h}\Delta_{A{},h}^{-1}\alpha =  \alpha. 
			\quad
			\mathrm{and} 
		\quad
			\Delta_{{A},h}^{-1}\Delta_{{A},h} \alpha =\alpha.
			$$
		\end{enumerate}
	\end{theorem}

 We stress that the same statement holds for $\Delta_{R,1}$ and $\Delta_{R,4}$, but not for $\Delta_{\G,h}$ (for any $h$) and not for $\Delta_{R,2}$ and $\Delta_{R,3}$, since in the latter cases the differential order is $12>Q$. However, since as stated in  Proposition \ref{hypo 2} they are homogeneous hypoelliptic operators, it is well-known that such 12-order operators admit an associated fundamental solution
 that belongs to $\mathcal{S}'(\G,E_0^h)$ of the form
 \begin{align}\label{fund solution }
     \widetilde{\mathcal{K}}=\mathcal{K}_0+p(x)\log(\rho(x))\,,
 \end{align}
 where  $\mathcal{K}_0$ is a vector-valued $\mc
			C^\infty(\G \setminus\{0\}$ function, homogeneous of degree $12-Q$, and $p(x)$ is a
homogeneous vector–valued polynomial of degree $12-Q$ (see e.g. Theorem 3.2.40-(b) in \cite{Fischer-Ruzhansky} and also Theorem 3.1 in \cite{VanSchaftingenYung}).

\medskip

In the sequel, we would need to estimate some higher-order derivatives of the fundamental solutions.
The following result is contained in Theorem 3.3. of \cite{VanSchaftingenYung}, and it ensures that the horizontal derivatives of order $\ell$ of the fundamental solution of $\Delta_{\G,h}$ can be seen associated to vector-valued kernels of type $12-\ell$.

\begin{theorem}[Theorem 3.3 in \cite{VanSchaftingenYung}]\label{global solution 12}
    Given the self-adjoint hypoelliptic operators $$\Delta_{\G,h}\colon \mathcal{D}'(\G,\mathbb R^{N_h})\to\mathcal{D}'(\G,\mathbb R^{N_h})$$ with $h=0,\cdots,5$, there exists $\mathcal{K}\in \mathcal{D}'(\G,\mathbb R^{N_h})$ so that for every $\alpha\in\mathcal{D}(\G,\mathbb R^{N_h})$, if we set \begin{align*}
        	\Delta_{\G,h}^{-1} \alpha:=
				\big(    
				\sum_{j}\alpha_j\ast  K_{1j},\dots,
				\sum_{j}\alpha_j\ast  K_{N_hj}\big)\,,
    \end{align*}
    and, as a shorthand notation we also write
    $$\Delta_{\G,h}^{-1} \alpha=\alpha\ast\mc K\,,$$
    then it holds that for any $\alpha \in \mc D(\G, \R^{N_h})$
    $$\Delta_{\G,h}\Delta_{\G,h}^{-1}\alpha=\Delta_{\G,h}(\alpha\ast\mc K)=\alpha.$$
    
    Moreover, if $\vert I\vert=\ell$, there exists $\mathcal{K}_I$ a vector-valued distribution so that $$X^I\alpha=\Delta_{\G,h}(\alpha\ast \mathcal{K}_I)$$ for every $\alpha\in \mathcal{D}(\G,\mathbb R^{N_h})$. If $\ell>2$, $\mathcal{K}_I$ is a homogeneous vector-valued distribution of degree $2-\ell$, and $\mathcal{K}_I$ agrees with a vector-valued  $C^\infty$-function on $\G\setminus\lbrace 0\rbrace$, so if further $\ell<12$, then $\mathcal{K}_I$ is a kernel of type $12-\ell$. 
\end{theorem}
\begin{remark} Roughly speaking, up to using right-invariant differential operators acting on $\mathcal{K}$, one can think of $\mathcal{K}_I$ as being the same as $X^I\mathcal{K}$ (see the proof of Theorem 3.3 in \cite{VanSchaftingenYung} for details).
    
\end{remark}

\begin{remark}\label{risultato di van schaftingen per altro laplaciano}
Applying again Theorem 3.3. of \cite{VanSchaftingenYung}, we can also say that the horizontal derivatives of order $\ell$ of the fundamental solution of both $\Delta_{R,2}$ and $\Delta_{R,3}$ can be seen associated to vector-valued kernels of type $12-\ell$ when $2<\ell<12$.
 \end{remark}

\subsection{A proof of the hypoellipticity of $\Delta_{A,h}$}
We prove now the hypoellipticity of the Rumin Laplacians introduced in Definition \ref{rumin laplacian}. Moreover, we will also state the div-curl inequalities associated with the Rumin complex $(E_0^\ast,d_c)$ for the Cartan group, when considering these particular Laplacians.

\subsubsection{Preliminaries on irreducible group representations}

\begin{definition}[Group representations]\label{Group representations}
    A representation $\pi$ of a Lie group $G$ on a Hilbert space ${H}_\pi\neq \lbrace 0\rbrace$ is a homomorphism $\pi$ of $G$ into the group of bounded linear operators on ${H}_\pi$ with bounded inverse. This means that
\begin{itemize}
    \item for every $x\in G$, the linear mapping $\pi(x)\colon {H}_\pi\to{H}_\pi$ is bounded and has bounded inverse,
    \item for any $x,y\in G$, we have $\pi(xy)=\pi(x)\pi(y)$.
\end{itemize}

A representation $\pi$ is unitary if $\pi(x)$ if unitary for every $x\in G$, i.e. $\pi\colon G\to U({H}_\pi)$. 
A representation $\pi$ is called irreducible when it has no closed invariant subspace.

\end{definition}
\begin{definition} Let $G$ be a simply connected Lie group and let $\pi$ be a representation of $G$ on a
Hilbert space $H_\pi$. A vector $v\in H_\pi$ is said to be smooth or of type $C^\infty$ if the
function 
$$G\ni x\to \pi(x)v\in H_\pi$$
is of class $C^\infty$. 
We denote by $S_\pi$
 the space of all smooth vectors of $\pi$.
    
\end{definition} 
We refer to the book \cite{Fischer-Ruzhansky} for  more details on the subject of group representations.
\begin{remark}\label{ridotto}
The representation $\pi$ determines a representation $\pi$ of the Lie algebra $\mathfrak g$ of $G$ as linear maps $S_\pi\to S_\pi$ (which extends uniquely to the algebra  of all left-invariant differential operator of $G$). If $\pi$ is irreducible    then,  as pointed out in  Rockland (\cite{Rockland_tams78}), there is a unitary equivalence, taking $H_\pi$ as $L^2(\R^k)$ (for some integer $k$, possibly $0$) and 
$
	S_{\pi}= \mathcal S(\rn k). 
	$
   \end{remark}

In particular, our proof is concentrated on the laplacian $\Delta_{A,2}$.
Going back to Definition \ref{Group representations},
 since we have described the operators $\Delta_{A,h}$ as matrix-valued operators,  we can think of the operator $\pi(\Delta_{A,2}) $ as corresponding to a matrix of the same size whose entry are operators on $S_\pi$ obtained by applying $\pi$ to the corresponding entry
of $\Delta_{A,h}$.
	As in Remark \ref{ridotto}, we can then think $
	S_{\pi}= \mathcal S(\rn k) 
	$ for a suitable $k$. 
 The same argument applies to $\pi(d_c)$ and $\pi(\delta_c)$.

\begin{proposition}\label{hypo}

    The operators $\Delta_{A,h}$ are hypoelliptic for each $h=0,\ldots,5$.
\end{proposition}
\begin{proof}
    We use Rockland's approach as in \cite{HN}, \cite{glowacki} and \cite{christ_et_al}. The cases of $h\neq 2,3$ are already covered by \cite{DaveHallerBGG}, Lemma 2.14. Let us deal now with the case of $h=2$. The case of $h=3$ will follow directly from the fact that $\Delta_{A,3}=\star\Delta_{A,2}\star$.

 Let us consider $\pi$ a nontrivial irreducible representation of $\mathbb{G}$ on the Hilbert space $L^2(\R^k)$, which determines a representation $\pi$ of the Lie algebra $\mathfrak g$  on $S_\pi$. Then the hypoellipticity of the operator $\Delta_{A,2}$ is equivalent to
	the injectivity of $\pi\big( \Delta_{A,2}\big)$ on $S_\pi$ (here for sake of simplicity we denote by $S_\pi$ also $(\mathcal S(\rn k) )^3$. 
 
	Therefore, assume that $\pi (\Delta_{A,2})u=0$ with $u\in   (\mathcal S(\rn k))^3$, i.e.,
 \begin{equation}\label{lapl 0}
\pi(d_c)\pi(\delta_c)u+\pi(\delta_c)\pi(A_\Delta)\pi(d_c)u=0\,.
	\end{equation}
 
 Hence for any $u=(u_1,u_2,u_3)\in\big(\mathcal{S}(\mathbb R^k)\big)^3$, if we multiply \eqref{lapl 0} by $u$ and we integrate the identity on 
	$(\rn k)^{3}$, we have
 \begin{align*}
     0=&\int_\G\langle \pi (d_c)\pi(\delta_c)u,u\rangle dV+\int_\G\langle \pi(\delta_c)\pi(A_\Delta)\pi(d_c)u,u\rangle dV\\=&\int_\G\Vert \pi(\delta_c)u\Vert^2 dV+\int_\G\langle \pi(A_\Delta)\pi(d_c)u,\pi(d_c)u\rangle dV
 \end{align*}
	
	We consider the second addend and we integrate by parts. Since
	\begin{equation*}\begin{split}
			\pi(-\Delta_{0}) = \pi (-X_1^2) + \pi(-X_2^2)
			= \pi(X_1)^* \pi(X_1) +  \pi(X_2)^* \pi(X_2).
	\end{split}\end{equation*} We obtain
	\begin{equation*}\begin{split}
			\int_\G &
			\langle  
				\begin{pmatrix}
					\pi(-\Delta_{0}) & 0 &0  \\
					0 & \pi(-\Delta_{0})  &0\\ 0 & 0 & \pi(-\Delta_{0})
				\end{pmatrix}
				\pi(d_c)u,\pi(d_c){u}\rangle\, dV
			 \\&=\sum_{i=1}^2\int_\G\langle\pi(X_i)^*\pi(X_i) \cdot \pi(d_c)u,\pi(d_c)u\rangle dV=\sum_{i=1}^2\int_\G \Vert \pi(X_i)\pi(d_c)u\Vert^2 dV  
	\end{split}\end{equation*}
	
	Hence, \eqref{lapl 0} is equivalent to
 \begin{align}
     &\pi(X_i)\pi(d_c)u=0 \text{ for each } i=1,2 \text{ and } \label{riga 1 bis}\\ &\pi(\delta_c)u=0\,.\label{riga 2 bis}
 \end{align} Expressing \eqref{riga 1 bis} component wise $d_c\colon E_0^2\to E_0^3$ with $d_c=(L_{j,k})$, we get
 \begin{align}
     &\sum_{k=1}^3\pi(X_1)\pi(L_{j,k})u_k=0\ ,\ \forall\, j=1,2,3\,\text{ and }\label{riga 1}\\ &\sum_{k=1}^3\pi(X_2)\pi(L_{j,k})u_k=0\ ,\ \forall\, j=1,2,3\,.\label{riga 2}
 \end{align}
 Applying $\pi(X_1)^\ast$ to \eqref{riga 1}, $\pi(X_2)^\ast$ to \eqref{riga 2}, and taking the sum, we get
 \begin{align*}
     \pi(-\Delta_{0})\sum_{k=1}^3\pi(L_{j,k})u_k=0
 \end{align*}
 which implies that $\sum_{k=1}^3\pi(L_{j,k})u_k=0$ for all $j=1,2,3$ by the hypoellipticity of the subLaplacian $\Delta_{0}$, and shortly
 \begin{align}\label{d_cu=0}
     \pi(d_c)u=0\,.
 \end{align}
Hence, by \eqref{riga 2 bis} and \eqref{d_cu=0}, the injectivity  of $\pi(\Delta_{A, 2})$ on $S_\pi$ is proved if we show that,
for $u\in S_\pi$ belonging to $$\mathrm{Ker}(\pi(d_c))\cap \mathrm{Ker}(\pi(\delta_c)),$$ then $u=0$.

Since the operators $d_c$ and $\delta_c$ are all given in their explicit formulation in Section \ref{Rumin conti espliciti}, it is possible to show the injectivity by applying inductively the fact that $X_1$ and $X_2$ are bracket generating the whole Lie algebra (see e.g. the proof of Lemma 4.1.8 in \cite{Fischer-Ruzhansky}). 

However, in this discussion, we follow the more elegant approach as in the proof of Theorem 5.2 of \cite{rumin_grenoble}. Indeed, by \cite{rumin_grenoble}, proof of Theorem 5.2, there exists
$X\in\mathfrak g$ such that, for any $v\in  \big(\mathcal S(\rn k)\big)^{3}$,
\begin{equation}\label{lapl 4}
v=Q_X\pi(d_c)v + \pi(d_c)Q_Xv,
\end{equation}
where
$$
Q_X:=\pi(\Pi_{E_0}\Pi_E)P_X i_X\pi(\Pi_E\Pi_{E_0}).
$$
Here $P_X$ is the inverse of $\pi(\mc L_X)$, $\mc L_X$ being the Lie derivative
along $X$. The above identity says that the smooth cohomology of the complex $\pi(d_c) $ vanishes i.e.
\begin{equation}\label{null cohomology}
    \mathrm{Ker}(\pi(d_c))\cap S_\pi=\pi(d_c)(S_\pi)\,.
\end{equation}
Keeping in mind that $\mathcal{R}(\pi(d_c))^\perp=\ker(\pi(\delta_c))$, the equality \eqref{null cohomology} implies the injectivity of the system $\pi(d_c)+\pi(\delta_c)$ which, as shown above, is equivalent to prove the  injectivity of $\pi(\Delta_{A, 2})$, and we are done.

\end{proof}

\section{Div-curl estimates associated to the three possible laplacians}\label{section 5}

Given the previous results, we can prove the following limiting Sobolev-type estimates. The estimates presented in Theorem \ref{H2} relate to the fundamental solutions of the Laplacians $\Delta_{A,h}$ of Definition \ref{rumin laplacian}, while the ones in Theorem \ref{C2} relate to the properties of the fundamental solutions of the Laplacians $\Delta_{\G,h}$ and $\Delta_{R,h}$. In particular, the estimates we obtain in the two following Theorems differ in degrees $h=2,3$. Luckily, when considering differential forms that are either closed or co-closed, we obtain the same estimates using any Laplacian, as proved in Corollary \ref{H2cor}.

\begin{theorem}\label{H2} Denote by $(E_0^*,d_c)$ the complex of intrinsic forms in $\G$.
Then there exists $C>0$ such that for any $h$-form
$u\in \mc D(\G, E_0^h)$, $0\le h\le 5$, such that
$$
\left\{\begin{aligned}
 d_c u = f \\ \delta_c u = g&
\end{aligned}
\right.
$$
if we use Definition \ref{rumin laplacian}, we have
 \begin{align*}
&\|u\|_{L^{Q/(Q-1)}(\G)\hphantom{, E_0^5}}\le   C\|f \|_{L^1(\G, E_0^{1})} \quad & \mbox{if $h=0$;}
\\
&\|u\|_{L^{Q/(Q-3)}(\G, E_0^1)}\le  C\big( \|f \|_{L^1(\G,E_0^2)} + \|\delta_c d_cg \|_{\mc H^1(\G)}\big) \quad & \mbox{if $h=1$;}
\\
&\|u\|_{L^{Q/(Q-3)}(\G, E_0^2)}\le   C\big( \|\nabla_\G f \|_{L^1(\G, E_0^3)} + \|g \|_{L^1(\G, E_0^1)}\big) \quad& \mbox{if $h=2$;}\\
&\|u\|_{L^{Q/(Q-3)}(\G, E_0^3)}\le   C\big(\|f \|_{L^1(\G, E_0^4)} + \|\nabla_\G g \|_{L^1(\G, E_0^2)}  \big) \quad & \mbox{if $h=3$.}\\
&\|u\|_{L^{Q/(Q-3)}(\G, E_0^4)}\le   C\big( \|d_c\delta_cf \|_{\mc H^1(\G, E_0^5)} + \|g \|_{L^1(\G, E_0^3)}\big) \quad & \mbox{if $h=4$;}
\\
&\|u\|_{L^{Q/(Q-1)}(\G, E_0^5)}\le    C\|g \|_{L^1(\G,E_0^4)} \quad & \mbox{if $h=5$.}
\end{align*}
Here $\nabla_\G f$ and $\nabla_\G g$ denote the horizontal gradient applied component-wise.

\end{theorem} 

\begin{theorem}\label{C2} Denote by $(E_0^*,d_c)$ the complex of intrinsic forms in $\G$.
Then there exists $C>0$ such that for any $h$-form
$u\in \mc D(\G, E_0^h)$, $0\le h\le 5$, such that
$$
\left\{\begin{aligned}
 d_c u = f \\ \delta_c u = g&
\end{aligned}
\right.
$$
if we use Definition \ref{RS Laplacian} or the Definition \ref{DH Laplacian}, we have
 \begin{align*}
&\|u\|_{L^{Q/(Q-1)}(\G)\hphantom{, E_0^5}}\le   C\|f \|_{L^1(\G, E_0^{1})} \quad & \mbox{if $h=0$;}
\\
&\|u\|_{L^{Q/(Q-3)}(\G, E_0^1)}\le  C\big( \|f \|_{L^1(\G,E_0^2)} + \|\delta_c d_cg \|_{\mc H^1(\G)}\big) \quad & \mbox{if $h=1$;}
\\
&\|u\|_{L^{Q/(Q-6)}(\G, E_0^2)}\le C\big(\|d_c\delta_c f\|_{L^1(\G, E_0^3)} + 
\|d_cg\|_{L^1(\G, E_0^2)} \big) \quad& \mbox{if $h=2$;}\\
&\|u\|_{L^{Q/(Q-6)}(\G, E_0^3)}\le   C\big(\|\delta_c f \|_{L^1(\G, E_0^3)} + \|\delta_c d_cg \|_{L^1(\G, E_0^2)}  \big) \quad & \mbox{if $h=3$;}\\
&\|u\|_{L^{Q/(Q-3)}(\G, E_0^4)}\le   C\big( \|d_c\delta_cf \|_{\mc H^1(\G, E_0^5)} + \|g \|_{L^1(\G, E_0^3)}\big) \quad & \mbox{if $h=4$;}
\\
&\|u\|_{L^{Q/(Q-1)}(\G, E_0^5)}\le    C\|g \|_{L^1(\G,E_0^4)} \quad & \mbox{if $h=5$.}
\end{align*}

\end{theorem}

Since the proofs of both results require some lengthy algebraic steps, they will be included in Sections \ref{Proof them H2} and \ref{Proof them C2} respectively.

\medskip
\begin{remark}\label{caso scalare}
   It is worth explicitly noting that with each of the three definitions of Laplacian, when dealing with functions (i.e. 0-forms) we recover the well known Sobolev inequality (see \cite{FGaW}, \cite{CDG}, \cite{MSC}, \cite{GN}). In particular, this fact was not obvious when dealing with the 12-order Laplacian $\Delta_{\G, 0}$. The proof for this latter case is given in the next section.
\end{remark}

\medskip
We end this section by stating a last result that holds
when a differential form $u$ is closed or co-closed. Indeed, for closed or co-closed forms, some of the previous estimates can be sharpened as follows, according to the degree of the exterior differential $d_c$. Moreover, all the subsequent estimates coincide with each of the three definitions of Laplacian.
\begin{theorem}
    \label{H2cor} Denote by $(E_0^*,d_c)$ the complex of intrinsic forms in $\G$.
Then there exists $C>0$ such that for any $h$-form
$u\in \mc D(\G, E_0^h)$, $0\le h\le 5$, such that
$$
\left\{\begin{aligned}
 d_c u = f \\ \delta_c u = 0&
\end{aligned}
\right.
$$
using either Definition \ref{DH Laplacian}, Definition \ref{rumin laplacian} or Definition \ref{RS Laplacian}, we have always the same estimates (according to the differential order of $d_c^{(h)}$): 
 \begin{align*}
&\|u\|_{L^{Q/(Q-1)}(\G)\hphantom{, E_0^5}}\le   C\|f \|_{L^1(\G, E_0^{1})} \quad & \mbox{if $h=0$;}
\\
&\|u\|_{L^{Q/(Q-3)}(\G, E_0^1)}\le  C\|f \|_{L^1(\G,E_0^2)}  \quad & \mbox{if $h=1$;}
\\
&\|u\|_{L^{Q/(Q-2)}(\G, E_0^2)}\le   C \| f \|_{L^1(\G, E_0^3)}  \quad& \mbox{if $h=2$;}\\
&\|u\|_{L^{Q/(Q-3)}(\G, E_0^3)}\le   C\|f \|_{L^1(\G, E_0^4)} \quad & \mbox{if $h=3$.}\\
&\|u\|_{L^{Q/(Q-1)}(\G, E_0^4)}\le   C\|f \|_{L^1(\G, E_0^5)}  \quad & \mbox{if $h=4$;}
\end{align*}
Analogously, there exists $C>0$ such that for any $h$-form
$u\in \mc D(\G, E_0^h)$, $0\le h\le 5$, such that
$$
\left\{\begin{aligned}
 d_c u = 0 \\ \delta_c u = g&
\end{aligned}
\right.
$$
we have the following estimates (according to the differential order of $\delta_c^{(h)}$):
 \begin{align*}
&\|u\|_{L^{Q/(Q-1)}(\G, E_0^1)}\le  C \|g \|_{L^1(\G)} \quad & \mbox{if $h=1$;}
\\
&\|u\|_{L^{Q/(Q-3)}(\G, E_0^2)}\le   C \|g \|_{L^1(\G, E_0^1)} \quad& \mbox{if $h=2$;}\\
&\|u\|_{L^{Q/(Q-2)}(\G, E_0^3)}\le   C \| g \|_{L^1(\G, E_0^2)}  \quad & \mbox{if $h=3$.}\\
&\|u\|_{L^{Q/(Q-3)}(\G, E_0^4)}\le   C\|g \|_{L^1(\G, E_0^3)} \quad & \mbox{if $h=4$;}
\\
&\|u\|_{L^{Q/(Q-1)}(\G, E_0^5)}\le    C\|g \|_{L^1(\G,E_0^4)} \quad & \mbox{if $h=5$.}
\end{align*}

\end{theorem}

\subsection{Rephrasing the closedness of forms on the Cartan group}\label{rephrasing}

Let us quickly revise some notations about tensors on a general Carnot group $\G$, in order to state a result due to  Chanillo \& Van Schaftingen that we will need later on. 

\subsection{Tensors on $\G$} Given a Carnot group $\G$, the first layer $V_1$ of its Lie algebra $\mathfrak g$ will be here denoted as $\mathfrak g_1$. We are also assuming $m_1=\mathrm{dim}\,\mathfrak g_1$

For $k \in \N$, let $\otimes^k \mathfrak g_1$ be
the $k$-fold tensor product of $\mathfrak g_1$. We will write $\mathcal{I}_k$ for the index set $\{1,\cdots, m_1\}^k$, and
$$X^\otimes _I=X_{i_1}\otimes\cdots\otimes X_{i_k}\quad
 \mathrm{for}\ I = (i_1,\ldots, i_k)\in \mathcal{I}_k$$
so that $\{X^\otimes_I\}_{I\in \mathcal{I}_k}$ is a basis of $\otimes^k \mathfrak g_1$. Then we have  a linear surjection from $\otimes^k \mathfrak g_1$  to the
vector space of all homogeneous left-invariant linear partial differential operators of order
$k$ on $\G$ with real coefficients, given by $$X^\otimes_I\to X^I\,.$$ 

Moreover, we denote by  $\mathrm{Sym} (\otimes^k\mathfrak g_1 )$ the subspace of all symmetric tensors in $\otimes^k \mathfrak g_1$.
There is a symmetrization map $$\mathrm{Sym} : \otimes^k \mathfrak g_1 \to \mathrm{Sym} (\otimes^k\mathfrak g_1 )\,,$$ which is a linear surjection given by
$$X^\otimes_I\to \frac{1}{k!}\sum_{\sigma\in  \mathcal{S}_k}X^\otimes_{\sigma(I)}, \quad I\in \mathcal{I}_k$$
where here $\mathcal{S}_k$ denotes the symmetric group on $k$ elements, and $\sigma(I) := (i_{\sigma(1)},\cdots,i_{\sigma(k)} )$ if $I = (i_1,\cdots i_k)\in \mathcal{I}_k$ as above, and $\sigma\in \mathcal{S}_k$. 

Denoting  by $\mc D(\G ,\mathrm{Sym} (\otimes^k \mathfrak g_1 ))$ the subspace
of compactly supported smooth symmetric horizontal $k$-tensors,
 we recall now the following  result due to Chanillo \& Van Schaftingen that turns to be a
key step in our proof of Theorem \ref{H2}.

\begin{theorem}[\cite{CvS2009}, Theorem 5]\label{chanillo_van5} Let $k\ge 1$ and
 $$F\in 
L^1(\G, \otimes^k \mathfrak g_1 ), \quad \Phi \in \mc D(\G, \mathrm{Sym} (\otimes^k\mathfrak g_1 )).
$$

Suppose
that $F$ has vanishing generalised divergence, i.e.
$$
\sum_{i_1,\dots,i_k}  X_{i_k}\cdots X_{i_1} F_{i_1,\dots,i_k} = 0 \quad\mbox{in $\mc D'(\G )$}.
$$
Then
$$
\Big| \int_{\G}\scal{\Phi}{F}_{\otimes^k\mathfrak g_1}\, dp \Big| \le C 
\|F\|_{L^1(\G, \otimes^k \mathfrak g_1)} \|\nabla_{\G{}} \Phi\|_{L^Q(\G, \otimes^k \mathfrak g_1)}\,,
$$
where $\nabla_\G\Phi$ denotes the horizontal gradient applied component-wise.
\end{theorem}

Let $k\ge 1$ be fixed, and let $F\in 
L^1(\G, \otimes^k \mathfrak g_1 )$ belong to the space of horizontal $k$-tensors.
We can write
$$
F = \sum_{i_1,\dots,i_k} F_{i_1,\dots,i_k} X_{i_1}\otimes \cdots \otimes X_{i_k}.
$$
As mentioned before, $F$ can be identified with the differential operator
$$
u\to Fu := \sum_{i_1,\dots,i_k} F_{i_1,\dots,i_k} X_{i_1}\cdots X_{i_k} u.
$$

As in \cite{BF7} and \cite{BFP1}, our proof of Theorem \ref{H2} relies  on  the fact 
(precisely stated in Theorem \ref{pierre} below) that the components 
with respect to a given basis of closed forms in $E_0^h$
can be viewed as the components of a horizontal vector field  with vanishing ``generalized horizontal divergence''. More precisely, in our case where $\G$ is again the Cartan group, we have the have the following result.

\begin{proposition}
    \label{pierre} Let $\alpha=\sum_J \alpha_J\xi_J^h\in \mathcal{D}(\G,E_0^h)$,
 $1\le h\le 4$, be such that
$$
d_c\alpha = 0.
$$

Then each component of $\alpha$ is of the form $\alpha_J=C_JF_J$, where $C_J$ is a real constant and $F_J$ is the component of a horizontal symmetric $s$-tensor $F$ associated to a differential operator  \begin{align*}
    F=\sum_{J=1}^{\dim E_0^h}F_JX_{i_1}\cdots X_{i_s}
\end{align*}
of order $s$, where $s$ coincides with the order of the Rumin differential $d_c^h\colon E_0^h\to E_0^{h+1}$, i.e. $s=2$ if $h=2$ and $s=3$ otherwise.

Moreover,
 the differential operator $F$ has vanishing generalised divergence, that is it satisfies
\begin{equation}
 \label{div gen}   
\sum_{J=1}^{\dim E_0^h} X_{i_1}\cdots X_{i_s} F_J=0\,.
\end{equation}

\end{proposition}
In order to prove this result we will use the 
following classical Cartan's formula in $\G$ 
(see, e.g., \cite{helgason},  identity (9) p. 21, though with a different normalization
of the wedge product).

\begin{theorem}[Cartan's formula]\label{cartan} Let $\omega $ be a smooth $h$-form of $(\Omega^*,d)$ (the
usual de Rham's complex), and let
$Z_0,Z_1,\dots,Z_h$ be smooth vector fields in $\G$. Then
\begin{equation}\label{cartan eq}\begin{split}
&\Scal{d\omega}{Z_0\wedge\cdots\wedge Z_h}  =
\sum_{i=0}^{h} (-1)^{i} Z_i \Scal{\omega}{Z_0\wedge\cdots \hat Z_i\cdots
\wedge Z_h}
\\& \hphantom{xxx}+
\sum_{0\le i <j\le h} (-1)^{i+j}  \Scal{\omega}{\,[Z_i,Z_j]\wedge \cdots \wedge\hat Z_i\wedge\cdots\wedge \hat Z_j\cdots}.
\end{split}\end{equation}
\end{theorem}

\medskip 

\begin{remark}\label{utile a cartan}
    Let $\alpha=\sum_J \alpha_J\xi_J^h\in E_0^h$,
 $1\le h\le 4$, be such that
$$
d_c\alpha = 0.
$$
Then $d\Pi_E\alpha =0$, i.e. $\Pi_E\alpha$ is closed in the usual sense.
\end{remark}
\begin{proof}
   Since $0=d_c\alpha:= \Pi_{E_0}d\Pi_E\alpha= \Pi_{E_0}\Pi_Ed\alpha$, 
by Theorem \ref{main rumin new}, ii), if we apply $\Pi_E$
to this equation, we get 
\begin{equation}\label{closed}
0= \Pi_E  \Pi_{E_0}\Pi_Ed\alpha = \Pi_Ed\alpha = d \Pi_E \alpha,
\end{equation}
by Theorem \ref{main rumin new}, iv), i.e,
$\Pi_E\alpha$ is closed in the usual sense, as claimed.  
\end{proof}
\begin{remark}
    Notice that, given an arbitrary Rumin $h$-form $\alpha\in E_0^h$, the $h+1$-form $d\Pi_E\alpha$ is not in general a Rumin form, but it is still an element in $(\mc R (d_0))^\perp=E_0^\bullet\oplus (\ker d_0)^\perp$. 
\end{remark}

\begin{proof}[Proof of Proposition \ref{pierre}]

Keeping in mind the previous remarks, we are going to apply Cartan formula choosing a suitable simple $h$-vector such that its dual does not belong to $\mathcal{R}(d_0)$. We divide the proof according to the degree $h$ of the form.
   \begin{itemize}
   \item[case $h=1$:] Let us consider  $\alpha=\alpha_1\theta_1+\alpha_2\theta_2$ and $d_c\alpha=0$. By Remark \ref{utile a cartan}, we have $d\Pi_E\alpha=0$. We want to apply Theorem \ref{cartan} to $\omega=\Pi_E\alpha$. By \eqref{rumin 6}, we have
   \begin{align*}
        \Pi_E\alpha=&\alpha-d_0^{-1}d_1\alpha-d_0^{-1}d_2\alpha+(d_0^{-1}d_1)^2\alpha-d_0^{-1}d_3\alpha+d_0^{-1}d_2d_0^{-1}d_1\alpha+d_0^{-1}d_1d_0^{-1}d_2\alpha+\\&-(d_0^{-1}d_1)^3\alpha+\text{ terms of weight 5 and 6}=\alpha-d_0^{-1}d_1\alpha-d_0^{-1}d_2\alpha+(d_0^{-1}d_1)^2\alpha\,,
   \end{align*}
   since $\mathcal{R}(d_0^{(1)})=\Omega^{2,2}\oplus\Omega^{2,3}$ by Remark \ref{Im d0 sulle  forme} i), we have that $d_0^{-1}$ vanishes on $\Omega^{2,4}$, $\Omega^{2,5}$, and $\Omega^{2,6}$.  

We now use \eqref{cartan eq} applied to $d\omega=d\Pi_E\alpha=0$, taking $Z_0=X_4$, $Z_1=X_1$. 
Hence, keeping in mind that $[Z_0,Z_1]=[X_4,X_1]=0$, we can write
\begin{align*}
    &0=Z_0\Scal{\alpha}{Z_1}-Z_1\Scal{\alpha}{Z_0}-\big(Z_0\Scal{d_0^{-1}d_1\alpha}{Z_1}-Z_1\Scal{d_0^{-1}d_1\alpha}{Z_0}\big)+\\&-\big(Z_0\Scal{d_0^{-1}d_2\alpha}{Z_1}-Z_1\Scal{d_0^{-1}d_2\alpha}{Z_0}\big)+\big(Z_0\Scal{(d_0^{-1}d_1)^2\alpha}{Z_1}-Z_1\Scal{(d_0^{-1}d_1)^2\alpha}{Z_0}\big)\,.
\end{align*}
Arguing using the weights of the forms and the vector fields, we know that $w(d_0^{-1}d_1\alpha)=2$, while $w(Z_0)=3$ and $w(Z_1)=1$, and they are therefore orthogonal by Remark \ref{orthogonality}. Moreover, $w(d_0^{-1}d_2\alpha)=w((d_0^{-1}d_1)^2\alpha)=3$, and so they act trivially on $Z_1$. Finally, since $w(\alpha)=1$, the expression then simplifies to
\begin{align*}
    0=&Z_0\Scal{\alpha}{Z_1}+Z_1\Scal{d_0^{-1}d_2\alpha}{Z_0}-Z_1\Scal{(d_0^{-1}d_1)^2\alpha}{Z_0}\\=&X_4\alpha_1+X_1\Scal{d_0^{-1}(-X_3\alpha_1\theta_1\wedge\theta_3-X_3\alpha_2\theta_2\wedge\theta_3)}{X_4}-X_1\Scal{d_0^{-1}d_1d_0^{-1}(X_1\alpha_2-X_2\alpha_1)\theta_1\wedge\theta_2}{X_4}\\\underbrace{=}_{\text{by }\eqref{d0 cartan}}&X_4\alpha_1+X_1\Scal{X_3\alpha_1\theta_4+X_3\alpha_2\theta_5}{X_4}+X_1\Scal{d_0^{-1}d_1(X_1\alpha_2-X_2\alpha_1)\theta_3}{X_4}\\=&X_4\alpha_1+X_1X_3\alpha_1+X_1\Scal{d_0^{-1}\big(X_1(X_1\alpha_2-X_2\alpha_1)\theta_1\wedge\theta_3+X_2(X_1\alpha_2-X_2\alpha_1)\theta_2\wedge\theta_3\big)}{X_4}\\=&X_4\alpha_1+X_1X_3\alpha_1-X_1\Scal{X_1(X_1\alpha_2-X_2\alpha_1)\theta_4+X_2(X_1\alpha_2-X_2\alpha_1)\theta_5}{X_4}\\=&X_4\alpha_1+X_1X_3\alpha_1-X_1X_1(X_1\alpha_2-X_2\alpha_1)\,.
\end{align*}

Hence the differential operator
\begin{align*}
    F:=&(X_4+X_1X_3+X_1^2X_2)\alpha_1-X_1^3\alpha_2=(3X_1^2X_2+X_2X_1^2-3X_1X_2X_1)\alpha_1-X_1^3\alpha_2
\end{align*}
has vanishing divergence in the weak sense, by the formula above. Moreover, it is associated to the horizontal 3-tensor $F\in \mathcal{D}(\mathbb G,\otimes^3\mathfrak g_1)$
\begin{align}\label{3-tensore F}
    F:=&\frac{\alpha_1}{3}(X_1\otimes X_1\otimes X_2+X_1\otimes X_2\otimes X_1+X_2\otimes X_1\otimes X_1)-\alpha_2X_1\otimes X_1\otimes X_1\,.
\end{align}

In this case, $(Z_1\wedge Z_0)^\natural=\theta_1\wedge\theta_4$ belongs to $\Xi_0^2$ (see \eqref{base 2 forme di Rumin}), and so the action of $d\Pi_E$ on $X_4\wedge X_1$ could also be recovered from the explicit formulae of $d_c$ in Section \ref{Rumin conti espliciti}. However, we decided to include the computations also for $h=1$ as a simpler case, as compared to the more general case that we need to cover when considering the following step of $h=2$.

       \item[case $h=2$:]
   
Let us consider the case when $\alpha=\sum_{i=1}^3\alpha_i\xi_i^2$ and $d_c\alpha=0$, keeping in mind that
\begin{equation}\label{basis}\begin{split}
\xi_1^2 & = \theta_1\wedge\theta_4,\quad
\xi_2^2 = \frac{\theta_1\wedge\theta_5+\theta_2\wedge\theta_4}{\sqrt{2}},\quad
\xi_3^2 = \theta_2\wedge\theta_5
\ .
\end{split}\end{equation}
By Remark \ref{utile a cartan} we have $d\Pi_E\alpha=0$. Again, we apply Theorem \ref{cartan} to $\omega=\Pi_E\alpha$. 
By \eqref{rumin 6}, 
\begin{align*}
    \Pi_E\alpha=&\alpha-d_0^{-1}d_1\alpha-d_0^{-1}d_2\alpha+d_0^{-1}d_1d_0^{-1}d_1\alpha+ \text{ terms of weight 7 and 8}\\=&\alpha-d_0^{-1}d_1\alpha-d_0^{-1}d_2\alpha+d_0^{-1}d_1d_0^{-1}d_1\alpha\,.
\end{align*}
since $\mathcal{R}(d_0^{(2)})=span_{C^\infty(\G)}\lbrace\theta_1\wedge\theta_2\wedge\theta_3\rbrace\oplus \Omega^{3,5}\oplus\Omega^{3,6}$ by Remark \ref{Im d0 sulle  forme} ii), and so $d_0^{-1}$ vanishes on $\Omega^{3,7}$, and $\Omega^{3,8}$.

We use \eqref{cartan eq} applied to $d\omega=d\Pi_E\alpha=0$, taking
$Z_0=X_5$, $Z_1= X_1,  Z_{2} = X_3$. Notice that the corresponding 3-covector $\theta_1\wedge\theta_3\wedge\theta_5$ does not belong to $\Xi_0^3$, but it is in  $(\mathcal{R}(d_0))^\perp$. In this case, it is not sufficient to use the explicit computations of the $d_c$ contained in Subsection \ref{Rumin conti espliciti} to recover the action of $d\Pi_E$.

We can write
\begin{equation}\label{cartan 2bis}\begin{split}
& 0  =  
Z_0  \Scal{\alpha}{Z_1\wedge Z_2} - Z_1  \Scal{\alpha}{Z_0\wedge Z_2} + Z_2  \Scal{\alpha}{Z_0\wedge Z_1}
\\&
- \Big(Z_0  \Scal{d_0^{-1}d_1\alpha}{Z_1\wedge Z_2} - Z_1  \Scal{d_0^{-1}d_1\alpha}{Z_0\wedge Z_2} + Z_2  \Scal{d_0^{-1}d_1\alpha}{Z_0\wedge Z_1}\Big)\\&- \Big(Z_0  \Scal{d_0^{-1}d_2\alpha}{Z_1\wedge Z_2} - Z_1  \Scal{d_0^{-1}d_2\alpha}{Z_0\wedge Z_2} + Z_2  \Scal{d_0^{-1}d_2\alpha}{Z_0\wedge Z_1}\Big)\\&+ \Big(Z_0  \Scal{(d_0^{-1}d_1)^2\alpha}{Z_1\wedge Z_2} - Z_1  \Scal{(d_0^{-1}d_1)^2\alpha}{Z_0\wedge Z_2} + Z_2  \Scal{(d_0^{-1}d_1)^2\alpha}{Z_0\wedge Z_1}\Big)
\\&
- \Scal{\alpha}{[Z_0,Z_1]\wedge Z_2} + \Scal{\alpha}{[Z_0,Z_2]\wedge Z_1} - \Scal{\alpha}{[Z_1,Z_2]\wedge Z_0}
\\&
 -\Big( - \Scal{d_0^{-1}d_1\alpha}{[Z_0,Z_1]\wedge Z_2} + \Scal{d_0^{-1}d_1\alpha}{[Z_0,Z_2]\wedge Z_1} - \Scal{d_0^{-1}d_1\alpha}{[Z_1,Z_2]\wedge Z_0} \Big)\\&
 -\Big( - \Scal{d_0^{-1}d_2\alpha}{[Z_0,Z_1]\wedge Z_2} + \Scal{d_0^{-1}d_2\alpha}{[Z_0,Z_2]\wedge Z_1} - \Scal{d_0^{-1}d_2\alpha}{[Z_1,Z_2]\wedge Z_0} \Big)\\&
 +\Big( - \Scal{(d_0^{-1}d_1)^2\alpha}{[Z_0,Z_1]\wedge Z_2} + \Scal{(d_0^{-1}d_1)^2\alpha}{[Z_0,Z_2]\wedge Z_1} - \Scal{(d_0^{-1}d_1)^2\alpha}{[Z_1,Z_2]\wedge Z_0} \Big)\,.
\end{split}\end{equation}
This long expression simplifies greatly thanks to the particular choice of the vector fields $Z_0,Z_1,Z_2$. In fact $[Z_0,Z_1]=[Z_0,Z_2]=0$.

Moreover, arguing using the weights of the forms and the vector fields, we know that $w(\alpha)=4$, while $w(Z_0\wedge Z_2)=5$, and $w(Z_1\wedge Z_2)=3$, and so the expression simplifies to
 
\begin{equation}\label{cartan 2tris}\begin{split}
& 0  =  
 Z_2  \Scal{\alpha}{Z_0\wedge Z_1}
\\&
- \Big(Z_0  \Scal{d_0^{-1}d_1\alpha}{Z_1\wedge Z_2} - Z_1  \Scal{d_0^{-1}d_1\alpha}{Z_0\wedge Z_2} + Z_2  \Scal{d_0^{-1}d_1\alpha}{Z_0\wedge Z_1}\Big)\\&- \Big(Z_0  \Scal{d_0^{-1}d_2\alpha}{Z_1\wedge Z_2} - Z_1  \Scal{d_0^{-1}d_2\alpha}{Z_0\wedge Z_2} + Z_2  \Scal{d_0^{-1}d_2\alpha}{Z_0\wedge Z_1}\Big)\\&+ \Big(Z_0  \Scal{(d_0^{-1}d_1)^2\alpha}{Z_1\wedge Z_2} - Z_1  \Scal{(d_0^{-1}d_1)^2\alpha}{Z_0\wedge Z_2} + Z_2  \Scal{(d_0^{-1}d_1)^2\alpha}{Z_0\wedge Z_1}\Big)\\&
 - \Scal{\alpha}{[Z_1,Z_2]\wedge Z_0}
+ \Scal{d_0^{-1}d_1\alpha}{[Z_1,Z_2]\wedge Z_0} \Big)
 + \Scal{d_0^{-1}d_2\alpha}{[Z_1,Z_2]\wedge Z_0}\\& 
 - \Scal{(d_0^{-1}d_1)^2\alpha}{[Z_1,Z_2]\wedge Z_0}
\end{split}\end{equation}

Again, considering the weights of the forms $w(d_0^{-1}d_2\alpha)=w((d_0^{-1}d_1)^2\alpha)=6$, while the highest weight of a 2-vector considered here is 5. Finally, since $w(d_0^{-1}d_1\alpha)=5$ and $w(Z_0\wedge Z_1)=4$, the expression \eqref{cartan 2tris} further simplifies to
\begin{equation}\label{cartan 2quattro}\begin{split}
 0  = & 
 Z_2  \Scal{\alpha}{Z_0\wedge Z_1}
+Z_1  \Scal{d_0^{-1}d_1\alpha}{Z_0\wedge Z_2} \\&
 + \Scal{d_0^{-1}d_2\alpha}{[Z_1,Z_2]\wedge Z_0} 
 - \Scal{(d_0^{-1}d_1)^2\alpha}{[Z_1,Z_2]\wedge Z_0}
\end{split}\end{equation}

Since
$$d_0(\theta_4\wedge\theta_5)=-\theta_1\wedge\theta_3\wedge\theta_5+\theta_2\wedge\theta_3\wedge\theta_5$$
we obtain
\begin{align*}
d_0^{-1}d_1(d_0^{-1}d_1\alpha)=&\left(\frac{X_2^2}{{2}}\alpha_1+(-X_2X_1-X_1X_2)\frac{\alpha_2}{2\sqrt{2}}+X_1^2\frac{\alpha_3}{2}\alpha_2\right)\theta_4\wedge\theta_5\,,
\end{align*}
Finally, 
\begin{align*}
d_0^{-1}d_2\alpha=&d_0^{-1}\left(-{X_3}\alpha_1 \theta_1\wedge\theta_3\wedge\theta_4-\frac{X_3}{\sqrt{2}}\alpha_2(\theta_1\wedge\theta_3\wedge\theta_4+\theta_2\wedge\theta_3\wedge\theta_5)-X_3{\alpha_3}\theta_2\wedge\theta_3\wedge\theta_5\,\right)=0,
\end{align*}
Hence, \eqref{cartan 2quattro} became
\begin{align*}
0=&-X_3\frac{\alpha_2}{\sqrt{2}}+X_1^2\alpha_3-X_1X_2\frac{\alpha_2}{\sqrt{2}}
+\frac{X_2^2}{{2}}\alpha_1+(-X_2X_1-X_1X_2)\frac{\alpha_2}{2\sqrt{2}}+X_1^2\frac{\alpha_3}{2}\\=&X_1^2\frac{3\alpha_3}{2}+(-\frac{5}{2\sqrt{2}}X_1X_2+\frac{X_2X_1}{\sqrt{2}})\alpha_2+\frac{X_2^2}{2}\alpha_1\,,
\end{align*}
Therefore we can consider
the following differential operator 
$$F:=\frac{3\alpha_3}{2}X_1^2+\frac{\alpha_1}{2}X_2^2+\frac{\alpha_2}{\sqrt{2}}\left(-\frac{5}{2}X_1 X_2+X_2   X_1\right)$$
that satisfies \eqref{div gen}. Notice that
$F$ can be identified with
 the tensor $F\in\mathcal{D}(\G,\otimes^2\vetguno{1})$ as
\begin{equation}
 \label{G caso h=2 due}   
F:=\frac{3\alpha_3}{2}X_1\otimes X_1+\frac{\alpha_1}{2}X_2\otimes X_2+\frac{\alpha_2}{\sqrt{2}}\left(-\frac{5}{2}X_1\otimes X_2+X_2   \otimes X_1\right)\,.
\end{equation}
\item[case $h=3$] Let us consider the case when $\alpha=\sum_{i=1}^3\alpha_i\xi_i^3$ and $d_c\alpha=0$. As above, we have $d\Pi_E\alpha=0$ and by \eqref{rumin 6}
\begin{align*}
    \Pi_E\alpha=&\alpha-d_0^{-1}d_1\alpha-d_0^{-1}d_2\alpha-d_0^{-1}d_3\alpha+(d_0^{-1}d_1)^2\alpha+d_0^{-1}d_2d_0^{-1}d_1\alpha+\\&+d_0^{-1}d_1d_0^{-1}d_2\alpha-(d_0^{-1}d_1)^3\alpha=\alpha-d_0^{-1}d_1\alpha+(d_0^{-1}d_1)^2\alpha\,,
\end{align*}
since $\mathcal{R}(d_0^{(3)})= \Omega^{4,7}\oplus\Omega^{4,8}$ by Remark \ref{Im d0 sulle  forme} iii), and so $d_0^{-1}$vanishes on $\Omega^{4,9}$. 

Moreover, for any $i=1,2,3$ the fact that $w(\xi_i^3)=6$ implies that $\xi_i^3$, as a 3-covector, always contains $\theta_3$ (see \eqref{base 3 forme di Rumin}). As $\theta_3$ is the only covector of weight 2, the action of $d_2$ on such a $\alpha$ will necessarily be trivial.

We could use \eqref{cartan eq} applied to $d\Pi_E\alpha=0$, taking $Z_0=X_1$, $Z_1=X_3$, $Z_2=X_4$, $Z_3=X_5$. Instead, in this case we can argue straightforwardly using the explicit expression of $d_c$ acting on 3-forms as shown in Section \ref{Rumin conti espliciti}. Indeed, notice that the 4-covector $\theta_1\wedge\theta_3\wedge\theta_4\wedge\theta_5$ belongs to $\Xi_0^4$, and so the action of $d\Pi_E$ on $X_1\wedge X_3\wedge X_4\wedge X_5$ is recovered from the explicit formulae in Section \ref{Rumin conti espliciti}. 

Therefore, we get that the differential operator
\begin{align*}
    F:=&(3X_1X_2^2+X_2^2X_1-3X_2X_1X_2)\alpha_1+\sqrt{2}(-2X_1X_2X_1+X_2X_1^2)\alpha_2+X_1^3\alpha_3\\=&(X_1X_2^2+X_3X_2-X_5)\alpha_1+{\sqrt{2}}(-X_1^2X_2+X_4)\alpha_2+X_1^3\alpha_3
\end{align*}
which has vanishing divergence in the weak sense, has the corresponding horizontal 3-tensor $F\in \mathcal{D}(\mathbb G,\otimes^3\mathfrak g_1)$
\begin{align}\label{3-tensore F}
    F:=&\frac{\alpha_1}{3}(X_1\otimes X_2\otimes X_2+X_2\otimes X_1\otimes X_2+X_2\otimes X_2\otimes X_1)+\\&-\frac{\sqrt{2}\alpha_2}{3}(X_1\otimes X_1\otimes X_2+X_1\otimes X_2\otimes X_1+X_2\otimes X_1\otimes X_1)+\alpha_3X_1\otimes X_1\otimes X_1\,.
\end{align}
\item[case $h=4$] The case of 4-forms is much more simpler, given the expression of the differential $d_c\colon E_0^4\to E_0^5$, which is an operator of order 1 in the horizontal derivatives. Therefore, it suffices to take $F=-\alpha_1 X_2+\alpha_2X_1\in\mathcal{D}(\G,\mathfrak g_1)$ which has vanishing divergence in the usual sense. 
\end{itemize}
\end{proof}
\subsection{{Proof of Theorem \ref{H2}}}\label{Proof them H2}
We are now ready to prove one of our main theorems.

%
%

\begin{proof}[Proof of Theorem \ref{H2}] The case $h=0$ is well known (see \cite{FGaW,CDG,MSC,GN}), and hence by Hodge-star duality this also sets the case of $h=5$.
 Similarly, we can restrict our proof of the result to forms in $E_0^h$, with $h=1,2$, and obtain the cases $h=3,4$ by applying Hodge-star duality to the complex $(E_0^\ast,d_c)$. 
  
  \medskip
  
\noindent \textbf{Case}  $\mathbf{ h =1.}$ If $u, \phi \in \mc D(\G, E_0^1)$, we can write
\begin{equation}\label{representation 2}\begin{split}
\scal{u}{\phi}_{L^2(\G,E_0^1)} & = \scal{ u}{\Delta_{A,1}\Delta_{A,1}^{-1} \phi}_{L^2(\G,E_0^1)}
\\&
= \scal{ u}{\big(\delta_cd_c + (d_c\delta_c)^3\big)\Delta_{A,1}^{-1}\phi}_{L^2(\G,E_0^1)}.
\end{split}
\end{equation}
Consider now the first term in the previous sum, 
$$
 \scal{ u}{\delta_cd_c \Delta_{A,1}^{-1}\phi}_{L^2(\G,E_0^1)}=  \scal{d_c u}{d_c \Delta_{A,1}^{-1}\phi}_{L^2(\G,E_0^2)}.
$$
Moreover, since
 $f, d_c \Delta_{A,1}^{-1}\phi\in E_0^{2}$, we can write $f=\sum_{\ell=1}^3 f_\ell\xi_\ell^{2}\,$, 
 $\,d_c \Delta_{A,1}^{-1}\phi=\sum_{\ell=1}^3 (d_c \Delta_{A,1}^{-1}\phi)_\ell\xi_\ell^{2}$,
 and hence
 we can reduce ourselves to estimate
 \begin{equation}\label{toprove}
 \scal{f_\ell}{(d_c \Delta_{A,1}^{-1}\phi)_\ell}_{L^2(\G)}\quad\mbox{for $\ell=1,2,3$.}
 \end{equation}
 Consider now the 
   horizontal $2$-tensor $F\in \mc D (\G, \otimes^2 \mathfrak g_1)$ 
  as defined in \eqref{G caso h=2 due} $$ F:=\frac{3f_3}{2}X_1\otimes X_1+\frac{f_1}{2}X_2\otimes X_2+\frac{f_2}{\sqrt{2}}\left(-\frac{5}{2}X_1\otimes X_2+X_2   \otimes X_1\right)\,,$$
which, by Proposition \ref{pierre}, satisfies the hypotheses of Theorem \ref{chanillo_van5}.

Suppose $f_\ell = f_1$ and consider the horizontal 2-tensor $F$. 
We consider now the symmetric horizontal $2$-tensor
$\Phi$:
$$
\Phi:= (d_c \Delta_{A,1}^{-1} \phi)_1\big(X_2\otimes X_2\big),
$$
so that 
\begin{equation*}
\scal{f_1}{(d_c \Delta_{A,1}^{-1} \phi)_1}_{L^2(\G)} = 
\scal{F}{\Phi}_{L^2(\G, \otimes^2 \mathfrak g_1)}.
\end{equation*}
\medskip
By Theorem \ref{chanillo_van5} 
\begin{equation}\label{eq 1sir bis}
\begin{split}
\big| \scal{f_1}{(d_c \Delta_{A,1}^{-1} \phi)_1}_{L^2(\G)} | &
\le \|G\|_{L^1(\G,\otimes^2 \mathfrak g_1)}\|\nabla_{\G}(d_c \Delta_{A,1}^{-1}\phi)_1\|_{L^Q(\G)}
\\
&
\le \|f\|_{L^1(\G,E_0^2)}\|\nabla_{\G}d_c \Delta_{A,1}^{-1}\phi\|_{L^Q(\G,E_0^2)}.
\end{split}
\end{equation}
On the other hand, keeping in mind that $d_c$ has order 3, $\nabla_{\G}d_c \Delta_{A,1}^{-1}\phi$ can be expressed as a sum of terms 
with components of the form
$$
\phi_j\ast X^I\tilde K_{ij}\, \quad\mbox{with $d(I)=4$.}
$$
By Theorem \ref{global solution}, iv) and Proposition \ref{kernel}, $X^I\tilde K_{ij}$ are  kernels
of order 2, so that, by Theorem \ref{Folland 75 thm}, we have
\begin{equation}\label{eq 1bis}
|\scal{f_1}{(d_c \Delta_{A,1}^{-1}\phi)_1}_{L^2(\G)}| \le  C\|f \|_{L^1(\G,E_0^2)}\|\phi\|_{L^{Q/3}(\G,E_0^1)}.
\end{equation}
The same argument applies to
$f_2$ and $f_3$ using again the 2-tensor $F$, with a suitable choice for the symmetric tensor $\Phi$. 
We obtain eventually
\begin{equation}\label{eq 1ter}
|\scal{f}{d_c \Delta_{A,1}^{-1}\phi}_{L^2(\G,E_0^2)}|\le C\|f \|_{L^1(\G,E_0^2)}\|\phi\|_{L^{Q/3}(\G,E_0^1)}.
\end{equation}

Consider now the second term in \eqref{representation 2}
$$
\scal{u}{ (d_c\delta_c )^3\Delta_{A,1}^{-1} \phi}_{L^2(\G,E_0^1)}= \scal{ \delta_cd_c\delta_c u}{\delta_cd_c\delta_c \Delta_{A,1}^{-1}\phi}_{L^2(\G)}.
$$
By Theorem \ref{global solution}, formula \eqref{numero2},
keeping in mind that $\delta_c$  is an operator of order 1 in the horizontal derivatives
when acting on $E_0^1$,
the quantity $\delta_cd_c\delta_c \Delta_{A,1}^{-1}\phi$ can be written as a sum of terms such as
$$
\phi_j\ast X^I \tilde K_{ij}, \quad\mbox{with  $d(I)=3$.}
$$

On the other hand,
$$
\scal{ \delta_cd_c\delta_c u}{\phi_j\ast X^I \tilde K_{ij} }_{L^2(\G)} =
\scal{\delta_cd_cg}{\phi_j\ast X^I \tilde K_{ij} }_{L^2(\G)}
= \scal{\delta_cd_cg \ast \ccheck(X^I \tilde K_{ij}) }{\phi_j }_{L^2(\G)}\,.
$$
By H\"older,
$$
\vert \scal{\delta_cd_cg \ast \ccheck(X^I \tilde K_{ij}) }{\phi_j }_{L^2(\G)}\vert \le \|\delta_cd_cg \ast \ccheck(X^I \tilde K_{ij})\|_{L^{Q/Q-3}}\|\phi_j\|_{L^{Q/3}}.
$$
Notice the $X^I \tilde K_{ij}$'s  and hence the $\ccheck(X^I \tilde K_{ij})$'s are  kernels of type 3.
Thus, by Theorem 6.10 in \cite{folland_stein},
$$
|\scal{ \delta_c d_c\delta_c u}{\phi_j\ast X^I \tilde K_{ij} }_{L^2(\G)} |
\le C \|\delta_cd_cg\|_{\mc H^1(\G)} \|\phi\|_{L^{Q/3}(\G, E_0^1)}.
$$
Combining this estimate with the one in \eqref{eq 1ter}, we get eventually
$$
|\scal{u}{\phi}_{L^2(\G, E_0^1)} |  \le C\big(\|f\|_{L^1(\G , E_0^2)} + 
\|\delta_cd_cg\|_{\mc H^1(\G)} \big) \|\phi\|_{L^{Q/3}(\G, E_0^1)},
$$
and hence
$$
\|u\|_{L^{Q/(Q-3)}(\G, E_0^1)}\le C\big( \|f \|_{L^1(\G, E_0^2)} + \|\delta_cd_cg \|_{\mc H^1(\G)}\big)\,,
$$
which concludes the proof in the case $h=1$.

  \medskip
  
\noindent \textbf{Case}  $\mathbf{ h =2.}$ For sake of simplicity, we set $\mathcal{K}:=\Delta_{A,2}^{-1}$ throughout this part of the proof. If $u,\phi\in E_0^2$ are smooth compactly supported
forms, then we can write
\begin{equation}\label{representation 3}
\begin{split}
\scal{u}{\phi}_{L^2(\G,E_0^2)} & = \scal{ u}{\Delta_{A,2}\mc K \phi}_{L^2(\G,E_0^2)}
\\&
= \scal{ u}{(\delta_cA_\Delta d_c + d_c\delta_c)\mc K\phi}_{L^2(\G,E_0^2)}.
\end{split}
\end{equation}
Consider now the term
$$
 \scal{ u}{\delta_cA_\Delta d_c \mc K\phi}_{L^2(\G,E_0^2)}=  \scal{d_c u}{A_\Delta d_c \mc K\phi}_{L^2(\G,E_0^3)}.
$$
Let us write $f:= d_c  u$. We can write $f=\sum_{\ell=1}^3 f_\ell \xi^3_\ell$.
As above, $d_c f=0$, and again
$$
\scal{d_c u}{A_\Delta d_c \mc K\phi}_{L^2(\G,E_0^3)}=\sum_\ell \scal{f_\ell}{(A_\Delta d_c \mc K\phi)_\ell}_{L^2(\G)}.
$$

Since $d_c\mathcal{K}\phi=\sum_{\ell=1}^3(d_c\mathcal{K}\phi)_\ell\xi_\ell^3$ and $A_\Delta=\begin{pmatrix}
     \delta_cd_c& 0 & 0\\
     0 &\delta_cd_c&  0\\
     0 & 0 &\delta_cd_c
 \end{pmatrix}$, where the $d_c$ in the matrix entries are the differentials acting on functions (i.e. 0-forms) and have order 1.

Since $A_\Delta$ is a diagonal matrix, we have that $(A_\Delta d_c\mathcal{K}\phi)_\ell=\delta_c^{(1)}d_c^{(0)}(d_c^{(2)}\mathcal{K}\phi)_\ell$, where in the formula we are highlighting the degree of the forms the operators $\delta_c$ and $d_c$ are acting on. 

 Hence, componentwise, for each $\ell=1,2,3$, $$\scal{f_\ell}{(A_\Delta d_c^{(2)} \mc K\phi)_\ell}_{L^2(\G)}=\scal{f_\ell}{\delta_c^{(1)}d_c^{(0)}(d_c^{(2)}\mathcal{K}\phi)_\ell}_{L^2(\G)}=\scal{d_c^{(0)}f_\ell}{d_c^{(0)}(d_c^{(2)}\mathcal{K}\phi)_\ell}_{L^2(\G,E_0^1)}.$$
 As above, by Proposition \ref{pierre}, each $\alpha_j=(d_cf_\ell)_j=(X_1f_\ell\theta_1+X_2f_\ell\theta_2)_j=X_jf_\ell$ with $j=1,2$ is one of the components of a horizontal 3-tensor with vanishing horizontal generalised divergence, see \eqref{3-tensore F}.

 To achieve the estimate of
 $ \scal{(d_c^{(0)}f_\ell)_j}{(d_c^{(0)} (d_c^{(2)} \mc K\phi)_\ell)_j}_{L^2(\G)}$, we define a new horizontal 3-tensor $\Phi$ as before,
 so that 
 $$
 \scal{(d_cf_\ell)_j}{(d_c( d_c \mc K\phi)_\ell)_j}_{L^2(\G)} = 
 \scal{F}{\Phi}_{L^2(\G,\mathfrak \otimes^3g_1)}.
 $$
 By Theorem \ref{chanillo_van5}, denoting by $\nabla_\G f=(\nabla_{\G}f_1,\nabla_\G f_2,\nabla_\G f_3)$
 $$
 \big| \scal{(d_c^{(0)}f_\ell)_j}{(d_c^{(0)}( d_c^{(2)} \mc K\phi)_\ell)_j}_{L^2(\G)} | \le \|\nabla_\G f\|_{L^1(\G,E_0^3)}\|\nabla_{\G}(d_c^{(0)}( d_c^{(2)}\mc K\phi)_\ell)\|_{L^{Q}(\G,E_0^3)}
 $$
 On the other hand, $\nabla_{\G}\big( d_c^{(0)}( d_c^{(2)} \mc K\phi)_\ell\big)$ can be expressed as a sum of terms with components
 of the form
 $$
 \phi_j\ast X^I\tilde K_{ij}, \quad\mbox{with $d(I)=4$,}
 $$
 since $d_c^{(2)} :E_0^2\to E_0^3$ is an operator of order 2 and $\nabla_\G$ is of order 1 in the horizontal derivatives.
The same can be repeated for any $\ell=1,2,3$.

By Theorem \ref{global solution}, iv) and Proposition \ref{kernel}, $X^I\tilde K_{ij}$ are  kernels
 of type 2, so that, by Theorem \ref{Folland 75 thm}, we have
 \begin{equation*}
 |\scal{f_\ell}{ (A_\Delta d_c \mc K\phi)_\ell}_{L^2(\G)}  |  \le  C\|\nabla_\G f \|_{L^1(\G,E_0^3)}\|\phi\|_{L^{Q/3}(\G,E_0^2)}.
 \end{equation*}
 The same argument can be carried out for all the components of $f$, yielding
 \begin{equation}\label{eq 1ter 2}
 |\scal{f}{A_\Delta d_c \mc K\phi}_{L^2(\G,E_0^3)}| \le C\|\nabla_\G f \|_{L^1(\G,E_0^3)}\|\phi\|_{L^{Q/3}(\G,E_0^2)}.
 \end{equation}

Consider now the term 
$$\scal{u}{d_c\delta_c\mathcal{K}\phi}_{L^2(\G,E_0^2)}=\scal{\delta_cu}{\delta_c\mathcal{K}\phi}_{L^2(\G,E_0^1)}.$$
We have
\begin{equation*}\begin{split}
\scal{\delta_c  u}{\delta_c\mc K\phi}_{L^2(\G,E_0^1)}  =
\scal{\star\delta_c  u}{\star\delta_c\mc K\phi}_{L^2(\G,E_0^4)} 
= \scal{\star g}{\star\delta_c\mc K\phi}_{L^2(\G,E_0^4)} .
\end{split}\end{equation*}
We notice now that $\star g$ is a $d_c$-closed form in $E_0^4$. Indeed, $d_c\star g=\star \delta_c g=\star\delta_c^2u=0$ up to a sign. We can then apply Proposition \ref{pierre} to get
\begin{equation}\label{eq g}
\big|\scal{\star g}{\star\delta_c\mc K\phi}_{L^2(\G,E_0^4)} \big| \le \|g\|_{L^1(\G,E_0^1)}\|\nabla_{\G}\ast \delta_c\mc K\phi\|_{L^Q(\G,E_0^4)}
\end{equation}
As above, $\nabla_{\G}\star \delta_c \mc K\phi$ can be expressed as a sum of terms with components of the form
$$
\phi_j\ast X^I\tilde K_{ij}, \quad\mbox{with $d(I)=4$,}
$$
since $\delta_c :E_0^3\to E_0^2$ is an operator of order 3 in the horizontal derivatives.
By Theorem \ref{global solution} and Proposition \ref{kernel}, $X^I\tilde K_{ij}$ are  kernels
of type 2, so that, again by Theorem \ref{Folland 75 thm}, we have
\begin{equation*}
\big|\scal{ g}{\delta_c\mc K\phi}_{L^2(\G,E_0^1)} \big|
 \le  C\|g \|_{L^1(\G,E_0^1)}\|\phi\|_{L^{Q/3}(\G,E_0^2)}.
\end{equation*}

Combining this estimate with the one in \eqref{eq 1ter 2}, we get eventually
$$
| \scal{u}{\phi}_{L^2(\G, E_0^2)} | \le C\big(\|\nabla_\G f\|_{L^1(\G, E_0^3)} + 
\|g\|_{L^1(\G, E_0^1)} \big) \|\phi\|_{L^{Q/3}(\G, E_0^2)},
$$
and hence
$$
\|u\|_{L^{Q/(Q-3)}(\G, E_0^2)}\le C\big(\|\nabla_\G f\|_{L^1(\G, E_0^3)} + 
\|g\|_{L^1(\G, E_0^1)} \big).
$$

This achieves the proof of the theorem.

\end{proof}

\subsection{Proof of Theorem \ref{C2}} \label{Proof them C2} Let us now prove Theorem \ref{C2}.
\begin{proof}[Proof of Theorem \ref{C2} using the Laplacians $\Delta_{R,h}$] By definition, the two Laplacians $\Delta_{R,h}$ and $\Delta_{A,h}$ differ only on degrees $h=2,3$. The arguments we can use here are analogous to the ones used in the proof of the previous result, but we have differences coming from the fact that we are dealing with operators of order 12. Hence, for completeness, we provide the crucial steps involved in the case of $h=2$.

\noindent \textbf{Case}  $\mathbf{ h =2.}$ If $u,\phi\in E_0^2$ are smooth compactly supported
forms, then we can write
\begin{equation}\label{representation 3}
\begin{split}
\scal{u}{\phi}_{L^2(\G,E_0^2)} & = \scal{ u}{\Delta_{R,2}\Delta_{R,2}^{-1} \phi}_{L^2(\G,E_0^2)}
\\&
= \scal{ u}{[(\delta_cd_c)^3 + (d_c\delta_c)^2]\Delta_{R,2}^{-1}\phi}_{L^2(\G,E_0^2)}.
\end{split}
\end{equation}
Consider now the term
$$
 \scal{ u}{(\delta_c d_c)^3 \Delta_{R,2}^{-1}\phi}_{L^2(\G,E_0^2)}=  \scal{d_c\delta_cd_c u}{d_c\delta_c d_c \Delta_{R,2}^{-1}\phi}_{L^2(\G,E_0^3)}.
$$
Let us write $f:= d_c  u$, and hence $d_c\delta_cf$ is closed.

 With the same arguments as in the proof of Theorem \ref{H2}, by Proposition \ref{pierre} and Theorem \ref{chanillo_van5},
 $$
 \big| \scal{d_c\delta_cd_c u}{d_c\delta_cd_c  \Delta_{R,2}^{-1}\phi}_{L^2(\G,E_0^3)} | \le \| f\|_{L^1(\G,E_0^3)}\|\nabla_{\G} d_c\delta_c d_c \Delta_{R,2}^{-1}\phi\|_{L^{Q}(\G,E_0^3)}
 $$
 On the other hand, $\nabla_{\G} d_c\delta_c d_c \Delta_{R,2}^{-1}\phi$ can be expressed as 
 $$
 \phi\ast \mathcal{K}_I, \quad\mbox{with $d(I)=7$,}
 $$
 since $d_c^{(2)} :E_0^2\to E_0^3$ is an operator of order 2 in the horizontal derivatives.

 By Theorem \ref{global solution 12}, and Proposition \ref{kernel}, $\mc K_I$ is a  kernel
 of type 5, so that, by Theorem \ref{Folland 75 thm}, we have
 \begin{equation}\label{aperitivo c2}
 \big| \scal{d_c\delta_cf}{d_c\delta_cd_c  \Delta_{R,2}^{-1}\phi}_{L^2(\G,E_0^3)} |   \le  C\|d_c\delta_c f \|_{L^1(\G,E_0^3)}\|\phi\|_{L^{Q/6}(\G,E_0^2)}.
 \end{equation}

Consider now the term 
$$\scal{u}{(d_c\delta_c)^2 \Delta_{R,2}^{-1}\phi}_{L^2(\G,E_0^2)}=\scal{d_c\delta_cu}{d_c\delta_c \Delta_{R,2}^{-1}\phi}_{L^2(\G,E_0^2)}.$$
By imposing $\delta_cu=g$, we have that $d_cg$ is closed and we can again apply Proposition \ref{pierre} and Theorem \ref{chanillo_van5} to get
\begin{equation}\label{eq g C2}
\big|\scal{d_c g}{d_c\delta_c \Delta_{R,2}^{-1}\phi}_{L^2(\G,E_0^2)} \big| \le \|d_cg\|_{L^1(\G,E_0^2)}\|\nabla_{\G} d_c\delta_c \Delta_{R,2}^{-1}\phi\|_{L^Q(\G,E_0^2)}
\end{equation}
As above, $\nabla_{\G} d_c\delta_c \Delta_{R,2}^{-1}\phi$ can be expressed as 
$$
\phi\ast \mc K_I, \quad\mbox{with $d(I)=7$,}
$$
since $\delta_c :E_0^3\to E_0^2$ is an operator of order 3 in the horizontal derivatives.
By Theorem \ref{global solution 12} and Proposition \ref{kernel}, $\mathcal{K}_I$ is a  kernel
of type 5, so that, by Theorem \ref{Folland 75 thm}, we have
\begin{equation*}
\big|\scal{d_c g}{d_c\delta_c\Delta_{R,2}^{-1}\phi}_{L^2(\G,E_0^2)} \big|
 \le  C\|d_cg \|_{L^1(\G,E_0^2)}\|\phi\|_{L^{Q/6}(\G,E_0^2)}.
\end{equation*}

Combining this estimate with the one in \eqref{aperitivo c2}, we get eventually
$$
| \scal{u}{\phi}_{L^2(\G, E_0^2)} | \le C\big(\|d_c\delta_c f\|_{L^1(\G, E_0^3)} + 
\|d_cg\|_{L^1(\G, E_0^2)} \big) \|\phi\|_{L^{Q/6}(\G, E_0^2)},
$$
and hence
$$
\|u\|_{L^{Q/(Q-6)}(\G, E_0^2)}\le C\big(\|d_c\delta_c f\|_{L^1(\G, E_0^3)} + 
\|d_cg\|_{L^1(\G, E_0^2)} \big).
$$

This achieves the proof of the theorem.
    
\end{proof}

\begin{proof}[Proof of Theorem \ref{C2} using the Laplacians $\Delta_{\G,h}$]
The proof differs from the previous one only for $h=0,1$ and again uses the same arguments. Therefore, we only give a gist of the proof for $h=0$, since the proof of the case $h=1$ can be easily adapted from the corresponding one given right above.

For $u\in E_0^0$, $\delta_c u=0$ and $\Delta_{\G,0}=(\delta_c d_c)^6$. Keeping in mind that $d_c:E_0^0\to E_0^1$ is a differential operator of order $1$, given $u,\phi\in \mc D(\G)$ we start again from the identity
\begin{equation}\label{representation 6}
\begin{split}
\scal{u}{\phi}_{L^2(\G)} & = \scal{ u}{\Delta_{\G,0}\Delta_{\G,0}^{-1} \phi}_{L^2(\G)}
\\&
= \scal{ d_cu}{(\delta_c (\delta_cd_c)^5\Delta_{\G,0}^{-1}\phi}_{L^2(\G,E_0^1)}.
\end{split}
\end{equation}
Now, since $d_cu=f$ is closed, arguing as in the previous proofs by means of Proposition \ref{pierre}, Theorem \ref{chanillo_van5} we eventually get an estimate of the form
$$
| \scal{u}{\phi}_{L^2(\G)} | \le\| f\|_{L^1(\G, E_0^1)} \|\nabla_{\G} d_c(\delta_c d_c)^5\Delta_{\G,0}^{-1}\phi\|_{L^Q(\G,E_0^1)}
$$

\end{proof}

In the next subsection, we prove Theorem \ref{H2cor}, while in the subsequent Subsection \ref{ultima} we will present a way to obtain different results than the ones achieved in Theorem \ref{C2} by considering new $L^p$-type spaces.

\subsection{Proof of Theorem \ref{H2cor} and final remarks}\label{H2cor e final remarks}

Let us show now the main changes for the proof of Theorem \ref{H2cor}
\begin{proof}[Proof of Theorem \ref{H2cor} using Definition \ref{rumin laplacian}]
We show the result in the case of a co-closed form (i.e. $g=0$) only when $h=2$ and $h=4$ (since in the other degrees the result is an immediate consequence of Theorem \ref{H2} by imposing $g=0$).

Arguing as in proof of Theorem \ref{H2}, writing again the identities in \eqref{representation 3} we have only the term
$$
 \scal{ u}{\delta_cA_\Delta d_c \Delta_{A,2}^{-1}\phi}_{L^2(\G,E_0^2)}=  \scal{d_c u}{A_\Delta d_c \Delta_{A,2}^{-1}\phi}_{L^2(\G,E_0^3)}.
$$
We write $f= d_c  u$ hence
 $d_c f=0$, and again 
$$
\scal{d_c u}{A_\Delta d_c \Delta_{A,2}^{-1}\phi}_{L^2(\G,E_0^3)}
$$
The operator $A_\Delta d_c \Delta_{A,2}^{-1}$ is associated with kernels of type $2$, hence again by Theorem \ref{chanillo_van5}, being $\nabla_\G A_\Delta d_c \Delta_{A,2}^{-1}$ associated to a kernel of type 1
 \begin{equation}\label{eq 1ter 2 bis}
 |\scal{f}{A_\Delta d_c \Delta_{A,2}^{-1}\phi}_{L^2(\G,E_0^3)}| \le C\| f \|_{L^1(\G,E_0^3)}\|\phi\|_{L^{Q/2}(\G,E_0^2)}
 \end{equation}
and hence, by duality,  the inequality for the case $h=2$ is proved.

Let us consider now $h=4$. 
If $u $ is co-closed, starting from the same procedure,
we have to analyze the term
$$
\scal{u}{ (\delta_cd_c )^3\Delta_{A,4}^{-1}  \phi}_{L^2(\G,E_0^4)}= \scal{ d_c u}{d_c\delta_cd_c\delta_cd_c \Delta_{A,4}^{-1}\phi}_{L^2(\G,E_0^5)}.
$$
Keeping in mind that $d_c$  is an operator of order 1 in the horizontal derivatives when acting on $E_0^4$, roughly speaking
the quantity $d_c\delta_cd_c\delta_cd_c \Delta_{A,4}^{-1}$ is associated to a kernel of type $1$ and hence, after applying Theorem \ref{chanillo_van5} we get that $\nabla_\G d_c\delta_cd_c\delta_cd_c \Delta_{A,4}^{-1}$ is of type $0$.
Hence, the conclusion follows by the following inequality
$$|\scal{u}{ (\delta_cd_c )^3\Delta_{A,4}^{-1}  \phi}_{L^2(\G,E_0^4)}|\le \| f \|_{L^1(\G,E_0^5)}\|\phi\|_{L^{Q}(\G,E_0^4)}. $$
This concludes the proof in the case of co-closed forms $u\in E_0^h$.

When $u$ satisfies $d_cu=0$ and $\delta_c u=g$, we can argue in a similar way.

\end{proof}
We now prove again Theorem \ref{H2cor} using the second definition of the Laplacian operators for co-closed forms. Therefore, the only case that we are left to show is for forms of degree $h=2,3$, since the Laplacians $\Delta_A$ and $\Delta_R$ coincide when $h=4$. 
\begin{proof}[Proof of Theorem \ref{H2cor} using Definition \ref{RS Laplacian}]
    If $u,\phi\in E_0^2$ are smooth compactly supported
forms, then we can write
\begin{equation}\label{representation 3}
\begin{split}
\scal{u}{\phi}_{L^2(\G,E_0^2)} & = \scal{ u}{\Delta_{R,2}\Delta_{R,2}^{-1} \phi}_{L^2(\G,E_0^2)}
\\&
= \scal{ u}{[(\delta_cd_c)^3 + (d_c\delta_c)^2]\Delta_{R,2}^{-1}\phi}_{L^2(\G,E_0^2)}.
\end{split}
\end{equation}
Since we are assuming that $\delta_cu=0$, we need to estimate only the term
$$
 \scal{ u}{(\delta_c d_c)^3 \Delta_{R,2}^{-1}\phi}_{L^2(\G,E_0^2)}=  \scal{d_c u}{d_c(\delta_c d_c)^2 \Delta_{R,2}^{-1}\phi}_{L^2(\G,E_0^3)}.
$$

 With the same arguments as in the proof of Theorem \ref{C2}, by Proposition \ref{pierre} and Theorem \ref{chanillo_van5},
$$
 \big| \scal{d_cu}{d_c(\delta_cd_c)^2  \Delta_{R,2}^{-1}\phi}_{L^2(\G,E_0^3)} | \le \| f\|_{L^1(\G,E_0^3)}\|\nabla_{\G} d_c(\delta_c d_c)^2 \Delta_{R,2}^{-1}\phi\|_{L^{Q}(\G,E_0^3)}
$$
 On the other hand, $\nabla_{\G} d_c(\delta_c d_c)^2 \Delta_{R,2}^{-1}\phi$ can be expressed in the form
 $$
 \phi\ast \mathcal{K}_I, \quad\mbox{with $d(I)=11$,}
 $$
 since $d_c^{(2)} :E_0^2\to E_0^3$ is an operator of order 2 in the horizontal derivatives.

 By Theorem \ref{global solution 12}, and Proposition \ref{kernel}, $\mc K_I$ is a kernel
 of type 1, so that, by Theorem \ref{Folland 75 thm}, we have
 \begin{equation}
    \label{stima f h=2}
 |\scal{u}{ \phi}_{L^2(\G,E_0^2)}  |  \le  C\| f \|_{L^1(\G,E_0^3)}\|\phi\|_{L^{Q/2}(\G,E_0^2)}.
 \end{equation}
 Reasoning by duality, we obtain the claim.

 The case of degree $h=3$ is very similar and yields the following estimate
 \begin{align*}
     |\scal{u}{ \phi}_{L^2(\G,E_0^3)}  |  \le  C\| f \|_{L^1(\G,E_0^4)}\|\phi\|_{L^{Q/3}(\G,E_0^3)}\,,
 \end{align*}
 from which one obtains the required estimate by duality.

 The arguments to obtain the estimates for a closed form $u\in E_0^h$ are again very similar. 
\end{proof}

To avoid excessive repetitions, we omit the proof of Theorem \ref{H2cor} using the Laplacian $\Delta_{\G,4}$, since it is sufficient to follow the same idea.
\subsection{An alternative Gagliardo-Niremberg type estimate}\label{ultima}

It is possible to prove Gagliardo-Niremberg  estimates also considering some variants of $L^p$ type spaces, as already considered in \cite{BFTr}. Indeed, in Theorems \ref{H2} and \ref{C2}, the estimates obtained for $h=1,2$ (and by Hodge duality also when $h=3,4$) the right-hand-side of the estimates contains not only the terms $d_cu=f$ and $\delta_cu=g$, but also some derivatives of them due to homogeneity required by the Laplacians. In this section, our aim is instead to obtain Gagliardo-Niremberg type inequalities where on the right-hand side we only have the $L^1$-norm (or $\mathcal{H}^1$-norm) of $f$ and $g$ alone, i.e. without any of their derivatives (see the statement of Theorem \ref{H2 intersezione}). 

Let us now set the function spaces that we are going to consider for such inequalities. First,
if $p,q\in [1,\infty]$, we can define the Banach spaces
$$
L^{p,q}(\G) := L^p(\G)\cap L^q(\G)
$$
endowed with the norm
$$
\| u\|_{L^{p,q}(\G)} :=  (\|u\|_{L^p(\G)} ^2+ \|u\|_{L^q(\G)}^2)^{1/2},
$$
and $\mc D(\G)$ is dense in $L^{p,q}(\G)$.
Analogous spaces
of differential forms can be defined in the usual way.

 Using these $L^{p,q}$-spaces, we manage to obtain the following estimates by duality.

\begin{theorem}\label{H2 intersezione} Denote by $(E_0^*,d_c)$ the complex of intrinsic forms in $\G$.
Then there exists $C>0$ such that for any $h$-form
$u\in \mc D(\G, E_0^h)$, $1\le h\le 4$, such that
$$
\left\{\begin{aligned}
 d_c u = f \\ \delta_c u = g&
\end{aligned}
\right.
$$
we have
 \begin{align*}
&\|u\|_{L^{Q/(Q-3)}+L^{Q/(Q-1)}(\G, E_0^1)}\le  C\big( \|f \|_{L^1(\G,E_0^2)} + \|g \|_{\mc H^1(\G)}\big) \quad & \mbox{if $h=1$;}
\\
&\|u\|_{L^{Q/(Q-3)}+L^{Q/(Q-2)}(\G, E_0^2)}\le   C\big( \| f \|_{L^1(\G, E_0^3)} + \|g \|_{L^1(\G, E_0^1)}\big) \quad& \mbox{if $h=2$;}\\
&\|u\|_{L^{Q/(Q-3)}+L^{Q/(Q-2)}(\G, E_0^3)}\le   C\big(\|f \|_{L^1(\G, E_0^4)} + \| g \|_{L^1(\G, E_0^2)}  \big) \quad & \mbox{if $h=3$.}\\
&\|u\|_{L^{Q/(Q-3)}+L^{Q/(Q-1)}(\G, E_0^4)}\le   C\big( \|f \|_{\mc H^1(\G, E_0^5)} + \|g \|_{L^1(\G, E_0^3)}\big) \quad & \mbox{if $h=4$;}
\end{align*}

\end{theorem} 
In the previous statement, we have used the vector spaces $L^p(\G)+L^q(\G)$ that, roughly speaking,  can be identified with the duals of the spaces $L^{p,q}(\G)$. Indeed, we can consider the vector space $L^p(\G) + L^q(\G)$ 
  endowed  with the norm
\begin{equation*}\begin{split}
\| u  \|_{L^p(\G) + L^q(\G)} : = \inf
\{  & ( \|u_1\|_{L^p(\G)}^2 + \|u_2\|_{L^q(\G)}^2)^{1/2};  \; 
\\&
u_1\in L^p(\G), u_2 \in L^q(\G), u = u_1+u_2\},
\end{split}\end{equation*}
notice that $L^p(\G) + L^q(\G)\subset L^1_{\mathrm{loc}}(\G)$. 

Moreover, if $p,q\in (1,\infty)$
 and $p',q'$ are their conjugate exponents, $(L^{p,q}(\G))^*$ is isometrically equals to  $L^{p'}(\G) + L^{q'}(\G)$ (see e.g. \cite{librointerpolationofoperators}  exercice 6 pag. 175, and  \cite{folland_real}, exercises 3 and 4 p. 186, and also \cite{BFTr} Proposition 3.1).

\begin{proof}[Proof of Theorem \ref{H2 intersezione}]The arguments used for the proof of this result are analogous to the ones used several times in the previous theorems. To avoid repeating the same steps, we only give the gist of the proof for only one of the possible Laplacians considered in this paper, keeping in mind the final estimates hold for all of them. 

Case $h=1$. If $u,\phi\in E_0^1$ are smooth and compactly supported forms, we can write
\begin{align}
    \scal{u}{\phi}_{L^2(\G,E_0^1)}=\scal{u}{\left[\delta_cd_c+(d_c\delta_c)^3\right]\Delta_{R,1}^{-1}\phi}_{L^2(\G,E_0^1)}\,.
\end{align}
The estimate for the term 
$    \scal{d_cu}{d_c\Delta_{R,1}^{-1}\phi}_{L^2(\G,E_0^2)}
$ is already contained in \eqref{eq 1ter} and then we have also
\begin{equation}\label{eq 1ter nuova}
|\scal{f}{d_c \Delta_{R,1}^{-1}\phi}_{L^2(\G,E_0^2)}|\le C\|f \|_{L^1(\G,E_0^2)}\big(\|\phi\|_{L^{Q/3}(\G,E_0^1)}+\|\phi\|_{L^{Q}(\G,E_0^1)}).
\end{equation}

The estimate for the second addend is obtained by writing 
\begin{align*}
    \scal{u}{(d_c\delta_c)^3\Delta_{R,1}^{-1}\phi}_{L^2(\G,E_0^1)}=\scal{g}{\delta_c(d_c\delta_c)^2\Delta_{R,1}^{-1}\phi}_{L^2(\G)}\,.
\end{align*}
The quantity $\delta_c(d_c\delta_c)^2\Delta_{R,1}^{-1}\phi$ can be written as a sum of terms of the form
\begin{align*}
    \phi_j\ast X^I\tilde{{K}}_{ij}\ \text{ with }d(I)=5\,.
\end{align*}
On the other hand, $\scal{g}{\delta_c(d_c\delta_c)^2\Delta_{R,1}^{-1}\phi}_{L^2(\G)}=\scal{g\ast \ccheck(X^I\tilde{{K}}_{ij})}{\phi_j}_{L^2(\G)}$. Moreover, the terms $X^I\tilde K_{ij}$ are kernels of type 1. Hence, using the same argument of the proof of Theorem \ref{H2} in the case when $h=1$, we get
\begin{align*}
    \vert\scal{g}{\delta_c(d_c\delta_c)^2\Delta_{R,1}^{-1}\phi}_{L^2(\G)}\vert\le&\Vert g\Vert_{\mathcal{H}^1(\G)}\Vert \phi\Vert_{L^{Q}(\G,E_0^1)}\\
    \le&\Vert g\Vert_{\mathcal{H}^1(\G)}\left(\Vert \phi\Vert_{L^{Q}(\G,E_0^1)}+\Vert\phi\Vert_{L^{Q/3}(\G,E_0^1)} \right)\,.
\end{align*}
This estimate, together with the one for $f$ given by \eqref{eq 1ter nuova}, gives
$$
|\scal{u}{\phi}_{L^2(\G, E_0^1)} |  \le C\big(\|f\|_{L^1(\G , E_0^2)} + 
\|g\|_{\mc H^1(\G)} \big) \big(\|\phi\|_{L^{Q/3}(\G, E_0^1)}+\|\phi\|_{L^{Q}(\G, E_0^1)}),
$$
and hence, by duality,
$$
\|u\|_{L^{Q/(Q-3)}+L^{Q/(Q-1)}(\G, E_0^1)}\le C\big( \|f \|_{L^1(\G, E_0^2)} + \|g \|_{\mc H^1(\G)}\big)\,,
$$

     {Case}  ${ h =2.}$ If $u,\phi\in E_0^2$ are smooth compactly supported
forms, then we can write
\begin{equation}\label{representation 3}
\begin{split}
\scal{u}{\phi}_{L^2(\G,E_0^2)} & = \scal{ u}{\Delta_{R,2}\Delta_{R,2}^{-1} \phi}_{L^2(\G,E_0^2)}
\\&
= \scal{ u}{[(\delta_cd_c)^3 + (d_c\delta_c)^2]\Delta_{R,2}^{-1}\phi}_{L^2(\G,E_0^2)}.
\end{split}
\end{equation}
Consider now the term
$$
 \scal{ u}{(\delta_c d_c)^3 \Delta_{R,2}^{-1}\phi}_{L^2(\G,E_0^2)}=  \scal{d_c u}{d_c(\delta_c d_c)^2 \Delta_{R,2}^{-1}\phi}_{L^2(\G,E_0^3)}.
$$
that as already been estimated in \eqref{stima f h=2}.
Therefore we have also
 \begin{equation}\label{stima f finale}
 |\scal{f}{(\delta_c d_c)^3 \Delta_{R,2}^{-1}\phi}_{L^2(\G,E_0^3)}| \le C\| f \|_{L^1(\G,E_0^3)}\big(\|\phi\|_{L^{Q/2}(\G,E_0^2)}+\|\phi\|_{L^{Q/3}(\G,E_0^2)}\big).
 \end{equation}

Consider now the term 
$$\scal{u}{(d_c\delta_c)^2 \Delta_{R,2}^{-1}\phi}_{L^2(\G,E_0^2)}=\scal{\delta_cu}{\delta_cd_c\delta_c \Delta_{R,2}^{-1}\phi}_{L^2(\G,E_0^1)}.$$
We have
\begin{equation*}\begin{split}
\scal{\delta_c  u}{\delta_cd_c\delta_c \Delta_{R,2}^{-1}\phi}_{L^2(\G,E_0^1)}  
= \scal{\star g}{\star\delta_cd_c\delta_c \Delta_{R,2}^{-1}\phi}_{L^2(\G,E_0^4)} .
\end{split}\end{equation*}
As already noticed, $\star g$ is a $d_c$-closed form in $E_0^4$ and we can then apply Proposition \ref{pierre} and Theorem \ref{chanillo_van5} to get
\begin{equation}\label{eq g C2}
\big|\scal{ g}{\delta_cd_c\delta_c \Delta_{R,2}^{-1}\phi}_{L^2(\G,E_0^4)} \big| \le \|g\|_{L^1(\G,E_0^1)}\|\nabla_{\G}\ast \delta_cd_c\delta_c \Delta_{R,2}^{-1}\phi\|_{L^Q(\G,E_0^4)}
\end{equation}
As above, $\nabla_{\G}\star \delta_c d_c\delta_c \Delta_{R,2}^{-1}\phi$ can be expressed as 
$$
\phi\ast \mc K_I, \quad\mbox{with $d(I)=10$.}
$$

By Theorem \ref{global solution 12} and Proposition \ref{kernel}, $\mathcal{K}_I$ is a  kernel
of type 2, so that, by Theorem \ref{Folland 75 thm}, we have
\begin{align*}
\big|\scal{ g}{\delta_c\mc K\phi}_{L^2(\G,E_0^1)} \big|
 \le & C\|g \|_{L^1(\G,E_0^1)}\|\phi\|_{L^{Q/3}(\G,E_0^2)}\\ \le&  C\|g \|_{L^1(\G,E_0^1)}\big(\|\phi\|_{L^{Q/2}(\G,E_0^2)}+\|\phi\|_{L^{Q/3}(\G,E_0^2)}\big).
\end{align*}

Combining this estimate with the one in \eqref{stima f finale}, we get eventually
$$
| \scal{u}{\phi}_{L^2(\G, E_0^2)} | \le C\big(\| f\|_{L^1(\G, E_0^3)} + 
\|g\|_{L^1(\G, E_0^1)} \big) \|\phi\|_{L^{Q/3, Q/2}(\G, E_0^2)},
$$
and hence, by duality
$$
\|u\|_{L^{Q/(Q-3)}+L^{Q/(Q-2)}(\G, E_0^2)}\le C\big(\| f\|_{L^1(\G, E_0^3)} + 
\|g\|_{L^1(\G, E_0^1)} \big).
$$
By Hodge-star duality, we get the corresponding estimates for $h=3,4$.

This achieves the proof of the theorem.
\end{proof}

\section*{Acknowledgements} A.B. is supported by the University of Bologna, funds for selected research topics, PRIN2022  \textit{Regularity problems in sub-Riemannian structures}– CUP J53D23003760006 (ref. 2022F4F2LH), and by GNAMPA of INdAM (Istituto Nazionale di Alta Matematica “F.
Severi”), Italy. 

F.T. would like to thank the Centro di
Ricerca Matematica Ennio De Giorgi and the Scuola Normale Superiore for the hospitality and support.

\bibliographystyle{alpha.bst}

\bibliography{BTr_Cartan}
\end{document}